\documentclass[12pt]{amsart}
\usepackage{amsmath}
\usepackage{amssymb}
\usepackage{amsfonts}
\usepackage{amsthm}
\usepackage{amscd}
\usepackage{bm}
\newtheorem{thrm}{Theorem}
\newtheorem{coro}{Corollary}
\newtheorem{prop}{Proposition}
\newtheorem{thm}{Theorem}[section]

\newtheorem{pro}{Proposition}[section]
\newtheorem{cor}{Corollary}[section]

\newcommand{\Z}{\mathbf{Z}}
\newcommand{\Q}{\mathbf{Q}}

\begin{document}
\title{First cohomology groups of the automorphism group of a free group with coefficients in the
       abelianization of the IA-automorphism group}
\address{Department of Mathematics, Faculty of Science Division II, Tokyo University of Science,
         1-3, Kagurazaka, Shinjuku-ku, Tokyo 162-8601, Japan}
\email{takao@rs.tus.ac.jp}
%ams classifcation
\maketitle
\begin{center}
{\sc Takao Satoh} \\

\vspace{0.5em}

{\footnotesize Department of Mathematics, Faculty of Science Division II, Tokyo University of Science, \\
         1-3, Kagurazaka, Shinjuku-ku, Tokyo 162-8601, Japan}
\end{center}

\begin{abstract}
We compute an integral twisted first cohomology group of the automorphism group of a free group with coefficients in the
abelianization $V$ of the IA-automorphism group of a free group. In particular, we show that it is generated by two crossed homomorphisms
constructed with the Magnus representation and the Magnus expansion due to Morita and Kawazumi respectively.
As a corollary, we see that the first Johnson homomorphism does not extend to the automorphism group of a free group as a crossed homomorphism
if the rank of the free group is greater than $4$.
\end{abstract}
\section{Introduction}\label{Int}

Let $F_n$ be a free group of rank $n \geq 2$ with basis $x_1, \ldots , x_n$, and $\mathrm{Aut}\,F_n$ the automorphism group of $F_n$.
The study of the (co)homology groups of $\mathrm{Aut}\,F_n$ with trivial coefficients has been developed for the last twenty years by many authors.
There are several remarkable results.
Gersten \cite{Ger} showed $H_2(\mathrm{Aut}\,F_n,\Z)=\Z/2\Z$ for $n \geq 5$.
Hatcher and Vogtmann \cite{HaV} showed $H_q(\mathrm{Aut}\,F_n,\Q)=0$ for $n \geq 1$ and $1 \leq q \leq 6$, except for
$H_4(\mathrm{Aut} \,F_4, \Q)= \Q$.
Furthermore, recently Galatius \cite{Gal} showed that the stable integral homology groups of
$\mathrm{Aut}\,F_n$ are isomorphic to those of the
symmetric group $\mathfrak{S}_n$ of degree $n$.
In particular, from his results, we see that the stable rational homology groups $H_q(\mathrm{Aut}\,F_n,\Q)$ of $\mathrm{Aut}\,F_n$ are trivial
for $n \geq 2q+1$.

\vspace{0.5em}

In this paper, we consider twisted (co)homology groups of $\mathrm{Aut}\,F_n$ from a viewpoint of the study of the Johnson homomorphism
of $\mathrm{Aut}\,F_n$. Let $H$ be the abelianization of $F_n$. The group $\mathrm{Aut}\,F_n$ naturally acts on $H$ and its dual group $H^* := \mathrm{Hom}_{\Z}(H,\Z)$.
There are a few computation for the (co)homology groups of $\mathrm{Aut}\,F_n$ with coefficients in $H$ and $H^*$.
Hatcher and Wahl \cite{HaN} showed that the stable homology groups of $\mathrm{Aut}\,F_n$ with coefficients in $H$ are trivial using
the stability of the homology groups of the mapping class groups of certain 3-manifolds.
In our previous papers \cite{Sa1} and \cite{Sa2},
we studied the stable twisted first and second (co)homology groups of $\mathrm{Aut}\,F_n$ with coefficients in $H$ and $H^*$, using the presentation for $\mathrm{Aut}\,F_n$
due to Gersten \cite{Ger}.
In particular, we obtained $H^1(\mathrm{Aut}\,F_n,H)=\Z$ for $n \geq 4$, and $H_2(\mathrm{Aut}\,F_n,H^*)=0$ for $n \geq 6$.

\vspace{0.5em}

Our research mentioned above is inspired by Morita's work for the mapping class group of a surface.
For $g \geq 1$, let $\Sigma_{g,1}$ be a compact oriented surface of genus $g$ with one boundary component, and $\mathcal{M}_{g,1}$ the mapping class
group of $\Sigma_{g,1}$. Namely, $\mathcal{M}_{g,1}$ is the group of isotopy classes of
orientation preserving diffeomorphisms of ${\Sigma}_{g,1}$ which fix the boundary pointwise.
The action of $\mathcal{M}_{g,1}$ on the fundamental group of ${\Sigma}_{g,1}$ induces a natural homomorphism
$\mathcal{M}_{g,1} \rightarrow \mathrm{Aut}\,F_{2g}$. It is known that this homomorphism is injective for any $g \geq 1$ due to classical works by Dehn and Nielsen.
Then we can consider $H$ as $\mathcal{M}_{g,1}$-modules for $n=2g$. We remark that $H^*$ is canonically isomorphic to $H$ by the Poincar\'{e} duality.
In \cite{Mo6}, Morita computed $H^1(\mathcal{M}_{g,1},H)=\Z$ for $g \geq 2$, and $H_2(\mathcal{M}_{g,1},H)=0$ for $g \geq 12$. (See also \cite{Mo0}.)
In particular, he showed that a crossed homomorphism induced from the Magnus representation of $\mathcal{M}_{g,1}$ generates $H^1(\mathcal{M}_{g,1},H)$.

\vspace{0.5em}

In general, the groups $\mathrm{Aut}\,F_n$ and $\mathcal{M}_{g,1}$ share many similar algebraic properties. 
If a certain result for either $\mathrm{Aut}\,F_n$ or $\mathcal{M}_{g,1}$ is obtained,
it would be natural to ask whether the corresponding result holds or not for the other.
As far as we compare the Morita's works with ours, it seems that $\mathrm{Aut}\,F_n$ and $\mathcal{M}_{g,1}$ behave similarly
with respect to the low dimensional twisted (co)homology groups.

\vspace{0.5em}

Here, we consider another $\mathrm{Aut}\,F_n$-module other than $H$ and $H^*$.
Let $\rho : \mathrm{Aut}\,F_n \rightarrow \mathrm{Aut}\,H$ be the natural homomorphism induced from the abelianization of $F_n$.
We identify $\mathrm{Aut}\,H$ with $\mathrm{GL}(n,\Z)$ by fixing a basis of $H$ induced from that of $F_n$.
The kernel of $\rho$ is called the IA-automorphism group of $F_n$, denoted by $\mathrm{IA}_n$.
The IA-automorphism group is a free group analogue of the Torelli subgroup of the mapping class group.
Although the study of the IA-automorphism group has a long history
since its finitely many generators were obtained by Magnus \cite{Mag} in 1935,
the combinatorial group structure of $\mathrm{IA}_n$ is still quite complicated.
For instance, any presentation for $\mathrm{IA}_n$ is not known in general.
Nielsen \cite{Ni0} showed that $\mathrm{IA}_2$ coincides with the inner automorphism group,
hence, is a free group of rank $2$.
For $n \geq 3$, however, $\mathrm{IA}_n$ is much larger than the inner automorphism group $\mathrm{Inn}\,F_n$.
Krsti\'{c} and McCool \cite{Krs} showed that $\mathrm{IA}_3$ is not finitely presentable.
For $n \geq 4$, it is not known whether $\mathrm{IA}_n$ is finitely presentable or not.
On the other hand, the abelianization $V$ of $\mathrm{IA}_n$ is completely determined by recent independent works of
Cohen-Pakianathan \cite{Co1, Co2}, Farb \cite{Far} and Kawazumi \cite{Kaw}. From their results, we have
$V \cong H^* \otimes_{\Z} \Lambda^2 H$ as a $\mathrm{GL}(n,\Z)$-module.

\vspace{0.5em}

Let $L$ be a principal ideal domain without $2$-torsion. (For example, $L=\Z$.)
In this paper, we determine the stable first cohomology group of $\mathrm{Aut}\,F_n$ with coefficients in $V_L:= V \otimes_{\Z} L$.
Here the ring $L$ is regarded as a trivial $\mathrm{Aut}\,F_n$-module. Our main theorem is
\begin{thrm}\label{T-I-1}($=$ Theorem {\rmfamily \ref{T-coh}}.)
For $n \geq 5$, if $L$ a principal ideal domain without $2$-torsion,
\[ H^1(\mathrm{Aut}\,F_n, V_L) = L^{\oplus 2}. \] 
\end{thrm}
We also show that the generators of $H^1(\mathrm{Aut}\,F_n, V_L)$ are constructed by the Magnus representation and the Magnus expansion
due to Morita \cite{Mo1} and Kawazumi \cite{Kaw} respectively. These are denoted by $f_M$ and $f_K$. (For details, see Section {\rmfamily \ref{S-Cro}}.)

\vspace{0.5em}

The computation of Theorem {\rmfamily \ref{T-I-1}} is motivated by a result for the mapping class group $\mathcal{M}_{g,1}$
due to Morita. In \cite{Mo4}, he computed the first cohomology group of $\mathcal{M}_{g,1}$ with coefficients in $\Lambda^3 H$, the free part of
the abelianization of the Torelli subgroup $\mathcal{I}_{g,1}$ of $\mathcal{M}_{g,1}$. In particular, he showed
$H^1(\mathcal{M}_{g,1},\Lambda ^3 H)={\Z}^{\oplus 2}$ for $g \geq 3$.
Hence, we also see that the corresponding result of $\mathcal{M}_{g,1}$ holds for $\mathrm{Aut}\,F_n$ in this case.

\vspace{0.5em}

In order to show the theorem above, we use Nielsen's presentation for $\mathrm{Aut}\,F_n$.
One of advantages of the generators-and-relations calculation is that 
by this method, we can determine $H^1(\mathrm{Aut}\,F_n, V_L)$ for many $L$ at the same time.
For example, $L=\Z$, $\Z/p\Z$ for any integer $p \in \Z$ such that $(p,2)=1$.

\vspace{0.5em}

Now, as an application of Theorem {\rmfamily \ref{T-I-1}}, we can see that the first Johnson homomorphism $\tau_1 : \mathrm{IA}_n \rightarrow V$ does not extend to
$\mathrm{Aut}\,F_n$ as a crossed homomorphism directly. (See \cite{Mo1} or \cite{Sat} for the definition of the Johnson homomorphism, for example.)
More precisely, for any principal ideal domain $L$, let denote $\tau_{1,L}$ the composition of the first Johnson homomorphism $\tau_1$ and the natural projection
$V \rightarrow V_L$. Then we have
\begin{coro}(= Corollary {\rmfamily \ref{p-joh}}.)
Let $L$ be a principal ideal domain without $2$-torsion. If $L$ does not contain $1/2$, then for $n \geq 5$,
there is no crossed homomorphism from $\mathrm{Aut}\,F_n$ to $V_L$ whose restriction to
$\mathrm{IA}_n$ coincides with $\tau_{1,L}$.
\end{coro}
We should remark that if a principal ideal domain $L$ contains $1/2$, then the first Johnson homomorphism
$\tau_1 : \mathrm{IA}_n \rightarrow V_L$ extends to a crossed homomorphism $\mathrm{Aut}\,F_n \rightarrow V_L$ due to Kawazumi \cite{Kaw}.
He explicitly construct a crossed homomorphism, denoted by $f_K$ in this paper,
which restriction to $\mathrm{IA}_n$ coincides with $\tau_{1,L}$ using the theory of Magnus expansions.
On the other hand, as to the mapping class group,
it has already known by Morita \cite{Mo5} that if $L$ contains $1/2$ then the first Johnson homomorphism
\[ \tau_1 : \mathcal{I}_{g,1} \rightarrow \Lambda^3 H \otimes_{\Z} L \]
of the mapping class group is uniquely extends to $\mathcal{M}_{g,1}$ as a crossed homomorphism where
$\mathcal{I}_{g,1}$ denotes the Torelli subgroup of $\mathcal{M}_{g,1}$.
Hence, we see that the groups $\mathrm{Aut}\,F_n$ and $\mathcal{M}_{g,1}$ also share a common property with respect to the extension of the first Johnson homomorphism.

\vspace{0.5em}

At the end of the paper, we consider the outer automorphism group $\mathrm{Out}\,F_n$.
In particular, we show

\begin{prop}(= Proposition {\rmfamily \ref{P-out}}.)
Let $L$ be a principal ideal domain without $2$-torsion. Then for $n \geq 5$,
\[ H^1(\mathrm{Out}\,F_n, V_L) = L. \]
\end{prop}

\vspace{0.5em}

This paper consists of six sections.
In Section {\rmfamily \ref{S-Pre}}, we fix some notation and conventions. Then we recall Nielsen's presentation for $\mathrm{Aut}\,F_n$.
In Section {\rmfamily \ref{S-Cro}}, we construct two crossed homomorphisms $f_M$ and $f_K$ from $\mathrm{Aut}\,F_n$ into $V_L$ for
any principal ideal domain $L$.
In Section {\rmfamily \ref{S-Coh}}, we compute the twisted first cohomology groups of $\mathrm{Aut}\,F_n$
using Nielsen's presentation.
In Section {\rmfamily \ref{S-App}}, we consider two applications. One is non-extendability of the Johnson homomorphism.
The other is a computation of the twisted first cohomology group of the outer automorphism group of a free group.

\tableofcontents

\section{Preliminaries}\label{S-Pre}

In this section, after fixing some notation and conventions,
we recall Nielsen's finite presentation for $\mathrm{Aut}\,F_{n}$, which is used to compute the
first cohomology groups in Section {\rmfamily \ref{S-Coh}}. Then we also recall
the IA-automorphism group of a free group and its abelianization.

\subsection{Notation and conventions}\label{Ss-Not}
\hspace*{\fill}\ 

\vspace{0.5em}

Throughout the paper, we use the following notation and conventions. Let $G$ be a group and $N$ a normal subgroup of $G$.
\begin{itemize}
\item The abelianization of $G$ is denoted by $G^{\mathrm{ab}}$.
\item The automorphism group $\mathrm{Aut}\,F_n$ of $F_n$ acts on $F_n$ from the right.
      For any $\sigma \in \mathrm{Aut}\,F_n$ and $x \in G$, the action of $\sigma$ on $x$ is denoted by $x^{\sigma}$.
\item For an element $g \in G$, we also denote the coset class of $g$ by $g \in G/N$ if there is no confusion.
\item Let $L$ be an arbitrary commutative ring. For any $\Z$-module $M$, we denote $M \otimes_{\Z} L$ by the symbol obtained by attaching a subscript $L$ to $M$,
      like $M_{L}$ or $M^{L}$. Similarly, for any $\Z$-linear map $f: A \rightarrow B$,
      the induced $L$-linear map $A_{L} \rightarrow B_{L}$ is denoted by $f_{L}$
      or $f^{L}$.
\item For elements $x$ and $y$ of $G$, the commutator bracket $[x,y]$ of $x$ and $y$
      is defined to be $[x,y]:=xyx^{-1}y^{-1}$.
\item For a group $G$ and a left $G$-module $M$, we set
\[\begin{split}
    \mathrm{Cros}(G,M) & := \{ f : G \rightarrow M \,|\, f : \mathrm{crossed} \,\,\, \mathrm{homomorphism} \}, \\
    \mathrm{Prin}(G,M) & := \{ g : G \rightarrow M \,|\, g : \mathrm{principal} \,\,\, \mathrm{crossed} \,\,\, \mathrm{homomorphism} \}.
  \end{split}\]
      Note that any modules of the crossed homomorphisms in this paper are considered to be left modules according to the usual custom.
      If $M$ is a right $\mathrm{Aut}\,F_n$-module induced from the action of automorphisms, we consider $M$ as a left module by the rule $\sigma \cdot m := m^{\sigma^{-1}}$
      for any $\sigma \in \mathrm{Aut}\,F_n$ and $m \in M$.
\end{itemize}

\subsection{Nielsen's Presentation}\label{Ss-Niel}
\hspace*{\fill}\ 

\vspace{0.5em}

For $n \geq 2$, let $F_n$ be a free group of rank $n$ with basis $x_1, \ldots , x_n$.
Let $P$, $Q$, $S$ and $U$ be automorphisms of $F_n$ given by specifying its images of the basis $x_1, \ldots, x_n$ as follows:

\vspace{0.5em}

\begin{center}
\begin{tabular}{c|c|c|c|c|c|c} \hline
           & $x_1$      & $x_2$ & $x_3$ & $\cdots$ & $x_{n-1}$ & $x_n$ \\ \hline
  $P$      & $x_2$      & $x_1$ & $x_3$ & $\cdots$ & $x_{n-1}$ & $x_n$ \\ 
  $Q$      & $x_2$      & $x_3$ & $x_4$ & $\cdots$ & $x_{n}$   & $x_1$ \\ 
  $S$      & $x_1^{-1}$ & $x_2$ & $x_3$ & $\cdots$ & $x_{n-1}$ & $x_n$ \\ 
  $U$      & $x_1 x_2$  & $x_2$ & $x_3$ & $\cdots$ & $x_{n-1}$ & $x_n$ \\ \hline
\end{tabular}
\end{center}

\vspace{0.5em}

In 1924, Nielsen \cite{Ni1} showed that the four elements above generate $\mathrm{Aut}\,F_{n}$.
Furthermore, he obtained the first finite presentation for $\mathrm{Aut}\,F_{n}$.
\begin{thm}[Nielsen \cite{Ni1}]
For $n \geq 2$, $\mathrm{Aut}\,F_{n}$ is generated by $P$, $Q$, $S$ and $U$
subject to finitely many relators. (For complete set of relators, see \cite{MKS} for example.)

\vspace{0,5em}

In this paper, in particular, we use the following relators: \\
\hspace{2em} {\bf{(N1)}}: $P^2$, $Q^n$, $S^2$, \\
\hspace{2em} {\bf{(N2)}}: $(QP)^{n-1}$, \\
\hspace{2em} {\bf{(N3)}}: $(PSPU)^2$, \\
\hspace{2em} {\bf{(N4)}}: $[S, QP]$, \\
\hspace{2em} {\bf{(N5)}}: $[U, Q^{-(l-1)} U Q^{l-1}]$ for $3 \leq l \leq n-1$, \\
\hspace{2em} {\bf{(N6)}}: $[U, Q^{-(l-1)} P Q^{l-1}]$ for $3 \leq l \leq n-1$, \\
\hspace{2em} {\bf{(N7)}}: $[P, Q^{-(l-1)} U Q^{l-1}]$ for $3 \leq l \leq n-1$, \\
\hspace{2em} {\bf{(N8)}}: $[U, Q^{-(n-2)} PUP^{-1} Q^{n-2}]$, \\
\hspace{2em} {\bf{(N9)}}: $U^{-1}PUP S U S P S $, \\
\hspace{2em} {\bf{(N10)}}: $(PQ^{-1}UQ)^2 UQ^{-1} U^{-1}QU^{-1}$. \\
It is easily cheched that the above elements are relators among the generators $P$, $Q$, $S$ and $U$.
\end{thm}

\vspace{0.5em}

Let $H$ be the abelianization of $F_n$, and $H^* :=\mathrm{Hom}_{\Z}(H,\Z)$ the dual group of $H$.
Let $e_1, \ldots , e_n$ be the basis of $H$ induced from $x_1, \ldots , x_n$, and $e_1^*, \ldots , e_n^*$
its dual basis of $H^*$. 
For any $\sigma \in  \mathrm{Aut}\,F_n$, the action of $\sigma$ on $e_i$ is given by
$\sigma \cdot e_i := x_i^{\sigma^{-1}} \in H$. In particular, the actions of $P$, $Q$, $S$ and $U$ on $e_i$ and $e_i^*$ are given by
\[\begin{tabular}{ll}
   $P \cdot e_k = \begin{cases}
                      e_2, \hspace{1em} k=1, \\
                      e_1, \hspace{1em} k=2, \\
                      e_k, \hspace{1em} k \neq 1, 2,
                    \end{cases}$ 
   & $P \cdot e_k^*  = \begin{cases}
                      e_2^*, \hspace{1em} k=1, \\
                      e_1^*, \hspace{1em} k=2, \\
                      e_k^*, \hspace{1em} k \neq 1, 2, 
                    \end{cases}$ \\
   $Q \cdot e_k = \begin{cases}
                      e_n, \hspace{1em} k=1, \\
                      e_{k-1}, \hspace{1em} k \neq 1,
                    \end{cases}$
   & $Q \cdot e_k^*  = \begin{cases}
                      e_n^*, \hspace{1em} k=1, \\
                      e_{k-1}^*, \hspace{1em} k \neq 1,
                    \end{cases}$ \\
   $S \cdot e_k = \begin{cases}
                      - e_1, \hspace{1em} k=1, \\
                      e_k, \hspace{1em} k \neq 1,
                    \end{cases}$
   & $S \cdot e_k^*  = \begin{cases}
                      - e_1^*, \hspace{1em} k=1, \\
                      e_k^*, \hspace{1em} k \neq 1,
                    \end{cases}$ \\
   $U \cdot e_k = \begin{cases}
                      e_1 - e_2, \hspace{1em} k=1, \\
                      e_k, \hspace{1em} k \neq 1,
                    \end{cases}$
   & $U \cdot e_k^*  = \begin{cases}
                      e_2^* + e_1^*, \hspace{1em} k=2, \\
                      e_k^*, \hspace{1em} k \neq 2.
                    \end{cases}$
  \end{tabular}\]

\subsection{IA-automorphism group}\label{Ss-IA}
\hspace*{\fill}\ 

\vspace{0.5em}

Here we recall the IA-automorphism group of a free group.
Fixing the basis $e_1, \ldots , e_n$ of $H$, we identify $\mathrm{Aut}\,H$ with $\mathrm{GL}(n,\Z)$.
The kernel of the natural homomorphism $\rho : \mathrm{Aut}\,F_n \rightarrow \mathrm{GL}(n,\Z)$ induced from
the abelianization of $F_n$ is called the IA-automorphism group of $F_n$, denoted by $\mathrm{IA}_n$.
Magnus \cite{Mag} showed that for any $n \geq 3$, the group $\mathrm{IA}_n$ is finitely generated by automorphisms
\[ K_{ij} : \begin{cases}
               x_i &\mapsto {x_j}^{-1} x_i x_j, \\
               x_t &\mapsto x_t, \hspace{4em} (t \neq i)
              \end{cases}\]
for distinct $i$, $j \in \{ 1, 2, \ldots , n \}$ and
\[  K_{ijk} : \begin{cases}
               x_i &\mapsto x_i x_j x_k {x_j}^{-1} {x_k}^{-1},  \\
               x_t &\mapsto x_t, \hspace{4em} (t \neq i)
              \end{cases}\] 
for distinct $i$, $j$, $k \in \{ 1, 2, \ldots , n \}$ such that $j>k$.

\vspace{0.5em}

Recently, Cohen-Pakianathan \cite{Co1, Co2}, Farb \cite{Far} and Kawazumi \cite{Kaw} independently determined
the abelianization of $\mathrm{IA}_n$. More precisely, they showed
\begin{equation}\label{CPFK}
\mathrm{IA}_n^{\mathrm{ab}} \cong H^* \otimes_{\Z} \Lambda^2 H
\end{equation}
as a $\mathrm{GL}(n,\Z)$-module. This abelianization is induced from the first Johnson homomorphism
\[ \tau_1 : \mathrm{IA}_n \rightarrow \mathrm{Hom}_{\Z}(H, \Lambda^2 H) = H^* \otimes_{\Z} \Lambda^2 H \]
defined by $\sigma \mapsto (x \mapsto x^{-1} x^{\sigma})$.
(For a basic material concerning the Johnson homomorphism, see \cite{Mo1} and \cite{Sat} for example.)
In this paper, we identify $\mathrm{IA}_n^{\mathrm{ab}}$ with $H^* \otimes_{\Z} \Lambda^2 H$ through $\tau_1$.
Then, (the coset classes of) the Magnus generators $K_{ij}$s and $K_{ijk}$s, whose images by $\tau_1$ are
\[ \tau_1(K_{ij}) = e_i^* \otimes e_i \wedge e_j, \hspace{1em} \tau_1(K_{ijk}) = e_i^* \otimes e_j \wedge e_k, \]
form a basis of $\mathrm{IA}_n^{\mathrm{ab}}$ as a free abelian group.
In the following, for simplicity, we write $V$ for $\mathrm{IA}_n^{\mathrm{ab}}$, and set
\[ \bm{e}_{j,k}^i := e_i^* \otimes e_j \wedge e_k \]
for any $i$, $j$ and $k$. Moreover, we consider the set
\[ I := \{ (i, j, k) \,|\, 1 \leq i \leq n, \hspace{1em} 1 \leq j < k \leq n \} \]
of the indices of the basis of $V$.

\vspace{0.5em}

Finally, we recall the inner automorphism group of $F_n$. A subgroup of $\mathrm{Aut}\,F_n$ consisting of automorphisms
given by conjugation by an element of $F_n$ is called the inner automorphism group of $F_n$. We denote it by $\mathrm{Inn}\,F_n$.
For each $1 \leq i \leq n$, set
\[ \iota_i := K_{1i} K_{2i} \stackrel{\check{i}}{\cdots} K_{ni} \in \mathrm{Inn}\,F_n. \]
It is easily seen that $\mathrm{Inn}\,F_n$ is a free group with basis $\iota_1, \ldots, \iota_n$.
We remark that the abelianization $(\mathrm{Inn}\,F_n)^{\mathrm{ab}}$ is naturally isomorphic to $H$ as a $\mathrm{GL}(n,\Z)$-module.
Furthermore, the inclusion $\mathrm{Inn}\,F_n \hookrightarrow \mathrm{IA}_n$ induces a
$\mathrm{GL}(n,\Z)$-equivariant injective homomorphism
\[ H = (\mathrm{Inn}\,F_n)^{\mathrm{ab}} \rightarrow \mathrm{IA}_n^{\mathrm{ab}} = H^* \otimes_{\Z} \Lambda^2 H \]
between their abelianizations.

\section{Construction of crossed homomorphisms}\label{S-Cro}

In this section, for any commutative ring $L$,
we introduce two crossed homomorphisms $f_M$ and $f_K$ from $\mathrm{Aut}\,F_n$ into $V_L= V \otimes_{\Z} L$,
due to Morita \cite{Mo1} and Kawazumi \cite{Kaw} respectively.
We remark that in their papers, the action of $\mathrm{Aut}\,F_n$ on $F_n$ is considered as the left one.
Hence, in this paper, whenever we use their notation and have to consider the left action of $\mathrm{Aut}\,F_n$ on $F_n$, we use
$\sigma(x) := x^{\sigma^{-1}}$ for any $\sigma \in \mathrm{Aut}\,F_n$ and $x \in F_n$.

\subsection{Morita's construction}
\hspace*{\fill}\ 

\vspace{0.5em}

First we construct a crossed homomorphism $f_M$ from $\mathrm{Aut}\,F_n$ into $V$ using the Magnus representation of $\mathrm{Aut}\,F_n$ due to
Morita \cite{Mo1}. Let
\[ \frac{\partial}{\partial x_j} : \Z[F_n] \longrightarrow \Z[F_n] \]
be the Fox's free derivations for $1 \leq j \leq n$. (For a basic material concerning with the Fox's derivation, see
\cite{Bir} for example.)
Let \ $\bar{ }: \Z[F_n] \rightarrow \Z[F_n]$ be the antiautomorphism induced from the map
$F_n \ni y \mapsto y^{-1} \in F_n$, and
$\frak{a} : \Z[F_n] \rightarrow \Z[H]$ the ring homomorphism induced from the abelianization $F_n \rightarrow H$.
For any matrix $A=(a_{ij}) \in \mathrm{GL}(n,\Z[F_n])$, set $A^{\mathfrak{a}} = (a_{ij}^{\mathfrak{a}}) \in \mathrm{GL}(n,\Z[H])$.
Then a map
\[ r_M : \mathrm{Aut}\,F_n \longrightarrow \mathrm{GL}(n,\Z[H]) \]
defined by
\[ \sigma \mapsto \biggl{(} \overline{\frac{\partial \sigma (x_j)}{\partial x_i}}
    {\biggl{)}}^{\frak{a}} \]
is called the Magnus representation of $\mathrm{Aut}\,F_n$.
We remark that $r_M$ is not a homomorphism but a crossed homomorphism. Namely, $r_M$ satisfies
\[ r_M(\sigma \tau) = r_M(\sigma) \cdot {r_M(\tau)}^{\sigma_*} \]
for any $\sigma$, $\tau \in \mathrm{Aut}\,F_n$ where ${r_M(\tau)}^{\sigma_*}$ 
denotes the matrix obtained from $r_M(\tau)$ by applying a ring homomorphism ${\sigma}_{*} : \Z[H] \rightarrow \Z[H]$
induced from $\sigma$ on each entry. (For detail for the Magnus representation, see \cite{Mo1}.)

\vspace{0.5em}

Now, observing the images of Nielsen's generators by $\mathrm{det} \circ r_M$, we verify that $\mathrm{Im}(\mathrm{det} \circ r_M)$ is contained in
a multiplicative abelian subgroup $\pm H$ of $\Z[H]$. In order to modify the image of $\mathrm{det} \circ r_M$, we consider the signature of $\mathrm{Aut}\,F_n$.
For any $\sigma \in \mathrm{Aut}\,F_n$, set $\mathrm{sgn}(\sigma) := \mathrm{det}(\rho(\sigma)) \in \{ \pm 1 \}$,
and define a map $f_M : \mathrm{Aut}\,F_n \longrightarrow \Z[H]$
by
\[ \sigma \mapsto \mathrm{sgn}(\sigma) \,\, \mathrm{det}(r_M(\sigma)). \]
Then the map $f_M$ is also crossed homomorphism which image of is contained in a multiplicative abelian subgroup $H$ in $\Z[H]$.
In the following, we identify the multiplicative abelian group structure of $H$ with the additive one.

\vspace{0.5em}

Finally, for any commutative ring $L$, by composing $f_M$ with a natural homomorphism
$H \rightarrow V \rightarrow V_L$ induced from the inclusion $\mathrm{Inn}\,F_n \hookrightarrow \mathrm{IA}_n$ and the projection $V \rightarrow V_L$,
we obtain an element in $\mathrm{Cros}(\mathrm{Aut}\,F_n, V_L)$, also denoted by $f_M$.

\subsection{Kawazumi's construction}
\hspace*{\fill}\ 

\vspace{0.5em}

Next, we construct another crossed homomorphism from $\mathrm{Aut}\,F_n$ into $V_L$ using the Magnus expansion of $F_n$ due to
Kawazumi \cite{Kaw}. (For a basic material for the Magnus expansion, see Chapter 2 in \cite{Bou}.)

\vspace{0.5em}

Let $\widehat{T}$ be the complete tensor algebra generated by $H$. For any Magnus expansion $\theta : F_n \rightarrow \widehat{T}$,
Kawazumi define a map
\[ \tau_1^{\theta} : \mathrm{Aut}\,F_n \rightarrow H^* \otimes_{\Z} H^{\otimes 2} \]
called the first Johnson map induced by the Magnus expansion $\theta$.
The map $\tau_1^{\theta}$ satisfies
\[ \tau_1^{\theta}(\sigma)([x]) = \theta_2(x) - |\sigma|^{\otimes 2} \theta_2(\sigma^{-1}(x)) \]
for any $x \in F_n$, where $[x]$ denotes the coset class of $x$ in $H$, $\theta_2(x)$ is the projection of $\theta(x)$ in $H^{\otimes 2}$,
and $|\sigma|^{\otimes 2}$ denotes the automorphism of $H^{\otimes 2}$ induced by $\sigma \in \mathrm{Aut}\,F_n$.
This shows that $\tau_1^{\theta}$ is a crossed homomorphism from $\mathrm{Aut}\,F_n$ to $H^* \otimes_{\Z} H^{\otimes 2}$.
In \cite{Kaw}, he also showed that $\tau_1^{\theta}$ does not depend on the choice of the Magnus expansion $\theta$, and that
the restriction of $\tau_1^{\theta}$ to $\mathrm{IA}_n$ is a homomorphism satisfying
\[ \tau_1^{\theta}(K_{ij}) = e_i^* \otimes e_i \otimes e_j - e_i^* \otimes e_j \otimes e_i,
   \hspace{1em} \tau_1^{\theta}(K_{ijk}) = e_i^* \otimes e_j \otimes e_k - e_i^* \otimes e_k \otimes e_j. \]

\vspace{0.5em}

Now, for any commutative ring $L$, compose $\tau_1^{\theta}$ and a natural projection
$H^* \otimes_{\Z} H^{\otimes 2} \rightarrow H^* \otimes_{\Z} \Lambda^2 H \rightarrow V_L$.
Then we obtain an element in $\mathrm{Cros}(\mathrm{Aut}\,F_n,V_L)$. In this paper, we denote it by $f_K$.
For $L=\Z$, from the result of Kawazumi as mentioned above, we see that the restriction of $f_K$ to $\mathrm{IA}_n$ coincides with
the double of the first Johnson homomorphism $\tau_1$. Namely, we have
\[ f_K(K_{ij}) = 2 \bm{e}_{i,j}^i, \hspace{1em} f_K(K_{ijk}) = 2 \bm{e}_{j,k}^i. \]
(See \cite{Mo1} or \cite{Sat} for the definition of the Johnson homomorphism, for example.)
Conversely, if $L$ contains $1/2$, the composition of the first Johnson homomorphism $\tau_1$ and the natural projection $V \rightarrow V_L$
extends to $\mathrm{Aut}\,F_n$ as a crossed homomorphism.

\subsection{Some observations}\label{Ss-som}
\hspace*{\fill}\ 

\vspace{0.5em}

In this subsection, we consider another crossed homomorphism $f_N$ in $\mathrm{Cros}(\mathrm{Aut}\,F_n, V_L)$ constructed from $f_M$ and $f_K$.
It is used to determine the first cohomology group $H^1(\mathrm{Aut}\,F_n,V_L)$ in Section {\rmfamily \ref{S-Coh}}.

\vspace{0.5em}

To begin with, we see the images of the crossed homomorphisms $f_M$ and $f_K$. From the definition, we have
\[ f_M(\sigma) := \begin{cases}
                    -(\bm{e}_{1,2}^2 + \bm{e}_{1,3}^3 + \cdots + \bm{e}_{1,n}^n), \hspace{1em} & \sigma=S, \\
                    0, \hspace{1em} & \sigma=P, Q, U
                 \end{cases}\]
and
\[ f_K(\sigma) :=\begin{cases}
                    -\bm{e}_{1,2}^1, \hspace{1em} & \sigma=U, \\
                    0, \hspace{1em} & \sigma=P, Q, S.
                 \end{cases}\]
These are obtained by straightforward calculations. We leave it to the reader as exercises.

\vspace{0.5em}

Next, for elements
\[ a_{j,k}^i := \begin{cases}
                  0, \hspace{1em} & i \neq j, \,\, k, \\
                  1, \hspace{1em} & i =k, \\ 
                  -1, \hspace{1em} & i =j,
                \end{cases}\]
in $L$ for $(i,j,k) \in I$, set
\[ \bm{a} := \sum_{(i,j,k) \in I} a_{j,k}^i \bm{e}_{j,k}^i \in V_L, \]
and let $f_{\bm{a}} \in \mathrm{Prin}(\mathrm{Aut}\,F_n,V_L)$ be a principal crossed homomorphism associated to $\bm{a} \in V_L$. 
Namely, for any $\sigma \in \mathrm{Aut}\,F_n$, it holds
\[\begin{split} f_{\bm{a}}(\sigma) & = \sigma \cdot \bm{a} - \bm{a}, \\
       & = \begin{cases}
             0, \hspace{1em} & \sigma=P, \, Q, \\
             -2 (\bm{e}_{1,2}^2 + \bm{e}_{1,3}^3 + \cdots + \bm{e}_{1,n}^n), \hspace{1em} & \sigma=S, \\
             \bm{e}_{1,2}^1 - (\bm{e}_{2,3}^3 + \bm{e}_{2,4}^4 + \cdots + \bm{e}_{2,n}^n), \hspace{1em} & \sigma=U. \\
           \end{cases} \\
  \end{split} \]
In fact, we have
\[\begin{split}
    f_{\bm{a}}(P) & = -(a_{1,2}^2+a_{1,2}^1) \bm{e}_{1,2}^1 + \sum_{k=3}^n (a_{2,k}^2 - a_{1,k}^1) \bm{e}_{1,k}^1 \\
                  & \hspace{1em} -(a_{1,2}^1+a_{1,2}^2) \bm{e}_{1,2}^2 + \sum_{k=3}^n (a_{1,k}^1 - a_{2,k}^2) \bm{e}_{2,k}^2 \\
                  & \hspace{1em} + \sum_{k=3}^n \Big{\{} (a_{2,k}^k -a_{1,k}^k) \bm{e}_{1,k}^k + (a_{1,k}^k -a_{2,k}^k) \bm{e}_{2,k}^k \Big{\}}, \\
                  & =0,
  \end{split}\]
\[\begin{split}
    f_{\bm{a}}(Q) & = \sum_{1 \leq j < i \leq n-1} (a_{j+1, i+1}^{i+1} - a_{j, i}^i) \bm{e}_{j,i}^i
                       + \sum_{j=1}^{n-1} (-a_{1, j+1}^1 - a_{j,n}^n) \bm{e}_{j,n}^n \\
                  & \hspace{1em} + \sum_{1 \leq i < j \leq n-1} (a_{i+1, j+1}^{i+1} - a_{i,j}^i) \bm{e}_{i,j}^i, \\
                  & = 0,
  \end{split}\]
\[\begin{split}
    f_{\bm{a}}(S) & = \sum_{2 \leq i \leq n} -2 a_{1, i}^i \bm{e}_{1,i}^i = -2 (\bm{e}_{1,2}^2 + \bm{e}_{1,3}^3 + \cdots + \bm{e}_{1,n}^n)
  \end{split}\]
and
\[\begin{split}
    f_{\bm{a}}(U) & = a_{1,2}^2 \bm{e}_{1,2}^1 + \sum_{3 \leq i \leq n} - a_{1, i}^i \bm{e}_{2,i}^i
                        + \sum_{3 \leq i \leq n} (a_{2, i}^{2} - a_{1,i}^1) \bm{e}_{2,i}^1, \\
                  & = \bm{e}_{1,2}^1 - (\bm{e}_{2,3}^3 + \bm{e}_{2,4}^4 + \cdots + \bm{e}_{2,n}^n).
  \end{split}\]

Now, we define
$f_N := 2 f_M - f_K - f_{\bm{a}} \in \mathrm{Cros}(\mathrm{Aut}\,F_n, V_L)$. From the arguments above, we have
\[ f_N(\sigma) =\begin{cases}
                    \bm{e}_{2,3}^3 + \bm{e}_{2,4}^4 + \cdots + \bm{e}_{2,n}^n, \hspace{1em} & \sigma=U, \\
                    0, \hspace{1em} & \sigma=P, Q, S.
                 \end{cases}\]
We use $f_N$ in Section {\rmfamily \ref{S-Coh}}.

\section{The first cohomology group}\label{S-Coh}

In the following, we always assume that $L$ is a principal ideal domain without $2$-torsion.
Set $V_L := V \otimes_{\Z} L$ as above. In this section, by using Nielsen's presentation for $\mathrm{Aut}\,F_n$, we show
\begin{thm}\label{T-coh}
For $n \geq 5$,
\[ H^1(\mathrm{Aut}\,F_n, V_L) = L^{\oplus 2}. \]
\end{thm}

\vspace{0.5em}

Here we give the outline of the computation. Let $F$ be a free group with basis $P$, $Q$, $S$ and $U$,
and $\varphi : F \rightarrow \mathrm{Aut}\,F_n$ the natural projection.
Then the kernel $R$ of $\varphi$ is a normal closure of the relators
{\bf (N1)}, $\ldots$, {\bf (N12)}.
Considering the five-term exact sequence of the Lyndon-Hochshild-Serre spectral sequence of the group extension
\[ 1 \rightarrow R \rightarrow F \rightarrow \mathrm{Aut}\,F_n \rightarrow 1, \]
we obtain an exact sequence
\[ 0 \rightarrow H^1(\mathrm{Aut}\,F_n, V_L) \rightarrow H^1(F, V_L) \rightarrow H^1(R, V_L)^F. \]
Observing this sequence at the cocycle level, we also obtain an exact sequence
\begin{equation}\label{eq-cross}
 0 \rightarrow \mathrm{Cros}(\mathrm{Aut}\,F_n,V_L) \rightarrow \mathrm{Cros}(F,V_L)
       \xrightarrow{\iota^*} \mathrm{Cros}(R,V_L)
\end{equation}
where $\iota^*$ is a map induced from the inclusion $\iota : R \hookrightarrow F$.
Hence we can consider $\mathrm{Cros}(\mathrm{Aut}\,F_n,V_L)$ as a subgroup consisting of elements of $\mathrm{Cros}(F,V_L)$
which are killed by $\iota^*$.
Hence, we can determine $\mathrm{Cros}(\mathrm{Aut}\,F_n,V_L)$ by using relators among the generators of Nielsen's presentation, and hence
$H^1(\mathrm{Aut}\,F_n, V_L)$.

\vspace{0.5em}

\textit{Proof of Theorem {\rmfamily \ref{T-coh}}.}
First, we consider the abelian group structure of $\mathrm{Cros}(F,V_L)$.
For any $\sigma \in F$ and a crossed homomorphism $f \in \mathrm{Cros}(F,V_L)$, we define elements $a_{j,k}^i(\sigma) \in L$ by
\[ f(\sigma) := \sum_{(i, j, k) \in I} a_{j,k}^i(\sigma) \bm{e}_{j,k}^i \in V_L. \]
Since $F$ is a free group generated by $P$, $Q$, $S$ and $U$, by the universality of a free group,
the crossed homomorphism $f$ is completely determined by
$a_{j,k}^i(\sigma)$ for $\sigma=P$, $Q$, $S$ and $U$.
More precisely, a map
\[ \mathrm{Cros}(F,V_L) \rightarrow L^{\oplus 2n^2(n-1)} \]
defined by
\[ f \mapsto \Big{(} a_{j,k}^i(P), \,\, a_{j,k}^i(Q), \,\, a_{j,k}^i(S), \,\, a_{j,k}^i(U) \Big{)}_{(i, j, k) \in I} \]
is an isomorphism between abelian groups. Through this map, we identify $\mathrm{Cros}(F,V_L)$ with $L^{\oplus 2n^2(n-1)}$ in the rest of the paper.

\vspace{0.5em}

In the following, we show that each $f \in \mathrm{Cros}(\mathrm{Aut}\,F_n,V_L) \subset \mathrm{Cros}(F,V_L)$ is determined by at most
\begin{equation}\begin{split}\label{eq-gen}
   & a_{j,k}^i(Q), \,\,\, 1 \leq i \leq n-1, \,\, 1 \leq j < k \leq n, \\
   & a_{j,k}^1(U), \,\,\, 1 \leq j < k \leq n, \\
   & a_{1,2}^2(S), \,\,\, a_{2,3}^3(U).
  \end{split}\end{equation}
Namely, define a map
\[ \Phi : \mathrm{Cros}(\mathrm{Aut}\,F_n, V_L) \rightarrow L^{\oplus (n^3-n^2+4)/2} \]
by
\[ f \mapsto \Big{(} (a_{j,k}^i(Q))_{i \neq n, \,\, 1 \leq j < k \leq n}, \,\, (a_{j,k}^1(U))_{1 \leq j < k \leq n}, \,\, a_{1,2}^2(S), \,\, a_{2,3}^3(U) \Big{)}. \]
Then our first goal is

\vspace{0.5em}

\noindent
{\bf Claim 1.} $\Phi$ is injective.

\vspace{0.5em}

Later, we see that elements in (\ref{eq-gen}) uniquely determine a crossed homomorphism from $\mathrm{Aut}\,F_n$ into $V_L$.
Our strategy to prove Claim 1 is as follows. Consider a commutative diagram of $L$-modules:
\[\begin{CD}
  0 @>>> \mathrm{Cros}(\mathrm{Aut}\,F_n,V_L) @>{\alpha}>> \mathrm{Cros}(F,V_L) @>{\iota^*}>> \mathrm{Cros}(R,V_L) \\
    @. @V{\Phi}VV @V{\cong}V{\beta}V @. \\
    @.  L^{\oplus (n^3-n^2+4)/2} @<{\gamma}<< L^{\oplus 2n^2(n-1)} @. \\
  \end{CD}\]
The first row is the exact sequaence (\ref{eq-cross}). The maps $\beta$, $\gamma$ are given by extraction of coefficients.
Then it suffices to show
$(\beta \circ \alpha)(\mathrm{Ker}(\Phi))=0$ in $L^{\oplus 2n^2(n-1)}$.
In other words, for any $f \in \mathrm{Cros}(\mathrm{Aut}\,F_n,V_L)$ such that all coefficients in (\ref{eq-gen}) are zero,
we show
\[ a_{j,k}^i(\sigma) = 0 \]
for any $\sigma= P$, $Q$, $S$ and $U$, and any $(i,j,k) \in I$. In order to do this, we use relators among the generators of Nielsen's presentation.
We divide the proof of Claim 1 into five steps. 

\vspace{0.5em}

\begin{flushleft}
{\bf Step I.} (Proof for $a_{j,k}^n(Q) = 0$.)
\end{flushleft}
From the relation {\bf (N1)}: $Q^n=1$, we obtain
\[ f(Q^n) = (1 + Q + Q^2 + \cdots + Q^{n-1}) f(Q) = 0. \]
For any $1 \leq j < k \leq n$, observing the coefficient of $\bm{e}_{j,k}^n$ in the equation above, we see
\[\begin{split}
   a_{j,k}^n(Q) & + a_{j+1,k+1}^1(Q) + a_{j+2,k+2}^2(Q) \cdots + a_{j+n-k,n}^{n-k}(Q) \\
     &  - a_{1,j+n-k+1}^{n-k+1}(Q) - \cdots - a_{k-j,n}^{n-j}(Q) + a_{1,k-j+1}^{n-j+1}(Q) + \cdots +a_{j-1,k-1}^{n-1}(Q) = 0,
  \end{split}\]
and hence $a_{j,k}^n(Q) = 0$. Therefore we see $a_{j,k}^i(Q) = 0$ for any $(i,j,k) \in I$.

\vspace{0.5em}

\begin{flushleft}
{\bf Step II.} (Some relations among $a_{j,k}^i(P)$ and $a_{j,k}^i(S)$.)
\end{flushleft}
Here we consider some linear relations among $a_{j,k}^i(P)$ and $a_{j,k}^i(S)$.

\vspace{0.5em}

From the relation {\bf (N1)}: $P^2=1$, we see
\[ f(P^2) = (1+P)f(P) = 0. \]
Observing the coefficients of $\bm{e}_{j,k}^{1}$, $\bm{e}_{1,k}^{1}$, $\bm{e}_{1,2}^{1}$, $\bm{e}_{2,k}^{1}$ and
$\bm{e}_{1,k}^{i}$ in the equation above, we obtain
\begin{gather}
 a_{j,k}^1(P) + a_{j,k}^2(P) = 0, \hspace{1em} 3 \leq j < k \leq n, \label{eq-p1} \\
 a_{1,k}^1(P) + a_{2,k}^2(P) = 0, \hspace{1em} 3 \leq k \leq n, \label{eq-p2} \\
 a_{1,2}^1(P) - a_{1,2}^2(P) = 0, \label{eq-p3} \\
 a_{2,k}^1(P) + a_{1,k}^2(P) = 0, \hspace{1em} 3 \leq k \leq n, \label{eq-p4} \\
 a_{1,k}^i(P) + a_{2,k}^i(P) = 0, \hspace{1em} 3 \leq i, k \leq n \label{eq-p5}
\end{gather}
respectively. 

\vspace{0.5em}

On the other hand, from the relation {\bf (N1)}: $S^2=1$, we see
\[ f(S^2) = (1+S)f(S) = 0. \]
Observing the coefficients of $\bm{e}_{j,k}^i$ for $2 \leq i \leq n$ and $2 \leq j < k \leq n$ in the equation above, we see
$2a_{j,k}^i(S) = 0$. Since $L$ does not contain any $2$-torsions, we obtain $a_{j,k}^i(S) = 0$.
Similarly, from the coefficients of $\bm{e}_{1,k}^1$ for $2 \leq k \leq n$, we see $a_{1,k}^1(S) = 0$.

\vspace{0.5em}

\begin{flushleft}
{\bf Step III.} (Proof for $a_{j,k}^i(U) = 0$.) This step consists of six parts.
\end{flushleft}

{\bf (i)} (Proof for $a_{1,2}^2(U) = a_{1,k}^i(U) = 0$ for $3 \leq i, k \leq n$.)
From the relation {\bf (N3)}: $(PSPU)^2=1$, we have
\begin{equation}\label{eq-r4}
\begin{split}
(PSPU + & 1)f(PSPU) \\
   & = (PSPU+1)(f(P) + P f(S) + PS f(P) + PSP f(U)) = 0
\end{split}
\end{equation}
The actions of $PS$, $PSP$ and $PSPU$ on $e_k$ and $e_k^*$ are given by
\[ PS \cdot e_k = \begin{cases}
                      - e_2, \hspace{0.5em} & k=1, \\
                      e_1, \hspace{0.5em} & k=2, \\
                      e_k, \hspace{0.5em} & k \neq 1, 2,
                    \end{cases} \hspace{1em}
   PS \cdot e_k^*  = \begin{cases}
                      - e_2^*, \hspace{0.5em} & k=1, \\
                      e_1^*, \hspace{0.5em} & k=2, \\
                      e_k^*, \hspace{0.5em} & k \neq 1, 2,
                    \end{cases} \]
\[ PSP \cdot e_k = \begin{cases}
                      -e_2, \hspace{0.5em} & k=2, \\
                      e_k, \hspace{0.5em} & k \neq 2,
                    \end{cases} \hspace{1em}
   PSP \cdot e_k^*  = \begin{cases}
                      - e_2^*, \hspace{0.5em} & k=2, \\
                      e_k^*, \hspace{0.5em} & k \neq 2,
                    \end{cases} \]
\[ PSPU \cdot e_k = \begin{cases}
                      e_1 + e_2, \hspace{0.5em} & k=1, \\
                      -e_2, \hspace{0.5em} & k=2, \\
                      e_k, \hspace{0.5em} & k \neq 1, 2,
                    \end{cases} \hspace{1em}
   PSPU \cdot e_k^*  = \begin{cases}
                      - e_2^* + e_1^*, \hspace{0.5em} & k=2, \\
                      e_k^*, \hspace{0.5em} & k \neq 2.
                    \end{cases} \]
Using this, we see that the coefficient of $\bm{e}_{2,k}^i$ for $3 \leq i, k \leq n$ in $(PSPU+1)f(PSPU)$
is equal to that of $\bm{e}_{1,k}^i$ in $f(PSPU)$, and to
\[ a_{1,k}^i(P) + a_{2,k}^i(S) + a_{2,k}^i(P) + a_{1,k}^i(U) = a_{1,k}^i(U) \]
by the argument above in Step II. Hence we obtain $a_{1,k}^i(U) = 0$.
Similarly, from the coefficient of $\bm{e}_{1,2}^1$ in $(PSPU+1)f(PSPU)$, which is equal to the coefficient of $\bm{e}_{1,2}^2$ in $f(PSPU)$ times $-1$, we see
\[ -(a_{1,2}^2(P) - a_{1,2}^1(S) - a_{1,2}^1(P) + a_{1,2}^2(U)) = - a_{1,2}^2(U) = 0. \]

\vspace{0.5em}

{\bf (ii)} (Proof for $a_{1,k}^2(U) = 0$ for $3 \leq k \leq n$.)
In general, if elements $\sigma, \tau \in \mathrm{Aut}\,F_n$ are commute, we have $f(\sigma) + \sigma f(\tau) = f(\tau) + \tau f(\sigma)$
from a relation $\sigma \tau=\tau \sigma$. Hence, we see
\begin{equation}
(\tau-1) f(\sigma) = (\sigma-1) f(\tau). \label{eq-comm}
\end{equation}
Consider the relator {\bf (N5)},
and applying (\ref{eq-comm}) for $\sigma=U$ and $\tau=Q^{-(l-1)} U Q^{l-1}$ for $3 \leq l \leq n-1$, we have
\[\begin{split}
    (Q^{-(l-1)} U & Q^{l-1} - 1) f(U) = (U-1) f(Q^{-(l-1)} U Q^{l-1}) \\
   & = (U-1)(f(Q^{-(l-1)}) + Q^{-(l-1)} f(U) + Q^{-(l-1)} U f(Q^{l-1})).
  \end{split}\]
Since $a_{j,k}^i(Q) = 0$ for any $(i,j,k) \in I$ by Step I, we see
\begin{equation}\label{eq-r2}
 (Q^{-(l-1)} U Q^{l-1} - 1) f(U) = (U-1) Q^{-(l-1)} f(U).
\end{equation}
Here the actions of $Q^{-(l-1)} U Q^{l-1}$ on $e_k$ and $e_k^*$ are given by
\[\begin{split} Q^{-(l-1)} U Q^{l-1} \cdot e_k & = \begin{cases}
                      e_l - e_{l+1}, \hspace{0.5em} & k=l, \\
                      e_k, \hspace{0.5em} & k \neq l,
                    \end{cases} \\
   Q^{-(l-1)} U Q^{l-1} \cdot e_k^* & = \begin{cases}
                      e_{l+1}^* + e_l^*, \hspace{0.5em} & k=l+1, \\
                      e_k^*, \hspace{0.5em} & k \neq l+1.
                    \end{cases}
  \end{split}\]
Observing this, we see the coefficient of $\bm{e}_{2,l}^{l+1}$ of $(Q^{-(l-1)} U Q^{l-1} - 1) f(U)$ is equal to $0$.
On the other hand, that of $(U-1) Q^{-(l-1)} f(U)$ is equal to
\begin{center}
 the coefficient of $\bm{e}_{1,l}^{l+1}$ of $Q^{-(l-1)} f(U)$ times $-1$,
\end{center}
and to
\begin{center}
 the coefficient of $\bm{e}_{1,n+2-l}^{2}$ of $f(U)$.
\end{center}
Hence, we see
\[ a_{1, n+2-l}^2(U) = 0 \]
for $3 \leq l \leq n-1$. Namely, $a_{1,k}^2(U) = 0$ for $3 \leq k \leq n-1$. 

\vspace{0.5em}

Similarly, consider the relator {\bf (N6)},
and applying (\ref{eq-comm}) for $\sigma=U$ and $\tau=Q^{-(l-1)} P Q^{l-1}$ for $3 \leq l \leq n-1$, we have
\begin{equation}\label{eq-r3}
 (Q^{-(l-1)} P Q^{l-1} - 1) f(U) = (U-1) Q^{-(l-1)} f(P).
\end{equation}
Here the actions of $Q^{-(l-1)} P Q^{l-1}$ on $e_k$ and $e_k^*$ are given by
\[\begin{split}
   Q^{-(l-1)} P Q^{l-1} \cdot e_k & = \begin{cases}
                      e_{l+1}, \hspace{0.5em} & k=l, \\
                      e_l, \hspace{0.5em} & k=l+1, \\
                      e_k, \hspace{0.5em} & k \neq l, l+1
                    \end{cases} \\
   Q^{-(l-1)} P Q^{l-1} \cdot e_k^* & = \begin{cases}
                      e_{l+1}^*, \hspace{0.5em} & k=l, \\
                      e_l^*, \hspace{0.5em} & k=l+1, \\
                      e_k^*, \hspace{0.5em} & k \neq l, l+1
                    \end{cases}
  \end{split}\]
Using this, from the coefficient of $\bm{e}_{1,l}^{2}$ in (\ref{eq-r3}), we see
\[ a_{1,l+1}^2(U)- a_{1,l}^2(U) = 0 \]
for $3 \leq l \leq n-1$. In particular, we have $a_{1,n}^2(U) = a_{1,n-1}^2(U) = 0$.

\vspace{0.5em}

{\bf (iii)} (Proof for $a_{j,k}^2(U) = 0$ for $2 \leq j < k \leq n$.)
First, we show $a_{2,k}^2(U) = 0$.
Now, observing the coefficients of $\bm{e}_{1,k}^{1}$ and $\bm{e}_{2,k}^{1}$ for $3 \leq k \leq n$ in $(PSPU + 1)f(PSPU)$, which are equal to
those of $2\bm{e}_{1,k}^1 + \bm{e}_{1,k}^2$ and $- \bm{e}_{2,k}^2 + \bm{e}_{1,k}^2$ in $f(PSPU)$ respectively, we obtain
\[\begin{split}
    2 a_{1,k}^1(P) + a_{1,k}^2(P) + 2 a_{2,k}^2(S) & + a_{2,k}^1(S) + 2 a_{2,k}^2(P) \\
         & - a_{2,k}^1(P) + 2a_{1,k}^1(U) - a_{1,k}^2(U) = 0, \\
    - a_{2,k}^2(P)  + a_{1,k}^2(P) - a_{1,k}^1(S)  & + a_{2,k}^1(S) - a_{1,k}^1(P) - a_{2,k}^1(P) \\
         & - a_{2,k}^2(U)  - a_{1,k}^2(U) = 0.
 \end{split}\]
Using (\ref{eq-p2}), and $a_{1,k}^1(S) = a_{2,k}^2(S) = 0$, we have
\begin{gather}
  a_{1,k}^2(P) + a_{2,k}^1(S)  - a_{2,k}^1(P) + 2a_{1,k}^1(U) - a_{1,k}^2(U) = 0, \label{eq-r5} \\
  a_{1,k}^2(P) + a_{2,k}^1(S) - a_{2,k}^1(P) - a_{2,k}^2(U)  - a_{1,k}^2(U) = 0. \label{eq-r6}
\end{gather}
Then, considering (\ref{eq-r5}) $-$ (\ref{eq-r6}), we obtain $a_{2,k}^2(U) = - 2a_{1,k}^1(U) = 0$ for $3 \leq k \leq n$.

\vspace{0.5em}

Next, we show $a_{j,k}^2(U) = 0$ for $3 \leq j < k \leq n$.
For $3 \leq j < l \leq n-1$, from the coefficient of $\bm{e}_{j,l+1}^2$ in (\ref{eq-r2}), we see
$a_{j,l}^2(U) = 0$.
Similarly, for $3 \leq l \leq n-2$, from the coefficient of $\bm{e}_{l+1,n}^2$ in (\ref{eq-r2}),
we see $a_{l,n}^2(U) = 0$.
To see $a_{n-1,n}^2(U) = 0$, we use (\ref{eq-r3}).
Observing the coefficient of $\bm{e}_{n-2,n}^2$ in (\ref{eq-r3}) for $l=n-2$,
we see $a_{n-1,n}^2(U) = a_{n-2,n}^2(U) = 0$.

\vspace{1em}

{\bf (iv)} (Proof for $a_{1,2}^i(U) = 0$ for $3 \leq i \leq n$.)
Observing the coefficient of $\bm{e}_{1,2}^l$ in (\ref{eq-r2}), we see
$a_{1,2}^{l+1}(U) = 0$ for $3 \leq l \leq n-1$.
To show $a_{1,2}^{3}(U) = 0$, considering the coefficient of $\bm{e}_{1,2}^3$ in (\ref{eq-r3}), we see
$a_{1,2}^{3}(U) = a_{1,2}^{4}(U) = 0$.

\vspace{1em}

{\bf (v)} (Proof for $a_{j,k}^i(U) = 0$ for $3 \leq i \leq n$ and $3 \leq j < k \leq n$.)
First, we consider the case where
$i \neq j, k$ and $j \neq k-1$, and show that $a_{j,k}^i(U) = a_{*,*}^n(U)$.
For $3 \leq l \leq n-1$ and $j, k \neq l, l+1$, from the coefficient of $\bm{e}_{j,k}^l$ in (\ref{eq-r3}), we see
\[ a_{j,k}^{l}(U) = a_{j,k}^{l+1}(U). \]
Similarly, observing the coefficients of $\bm{e}_{j,l+1}^l$ and $\bm{e}_{l+1,k}^l$ in (\ref{eq-r3}), we obtain
\begin{gather}
 a_{j,l+1}^l(U) = a_{j,l}^{l+1}(U), \hspace{1em} 3 \leq j < l \leq n-1, \label{eq-r13} \\
 a_{l+1,k}^l(U) = a_{l,k}^{l+1}(U), \hspace{1em} 3 \leq l < k-1 \leq n-1 \label{eq-r12}
\end{gather}
respectively. Using these equations, we obtain,
\begin{gather}
  a_{j,k}^{i}(U) = a_{j,k}^{i+1}(U) = \cdots = a_{j,k}^{n}(U), \hspace{1em} k+1 \leq i \label{eq-r9} \\
  a_{j,k}^{i}(U) = a_{j,k}^{i+1}(U) = \cdots = a_{j,k}^{k-1}(U) \stackrel{(\ref{eq-r13})}{=} a_{j,k-1}^{k}(U) \stackrel{(\ref{eq-r9})}{=}
    a_{j,k-1}^{n}(U), \hspace{1em} j+1 \leq i \leq k-1, \label{eq-r8} \\
  a_{j,k}^{i}(U) = a_{j,k}^{i+1}(U) = \cdots = a_{j,k}^{j-1}(U) \stackrel{(\ref{eq-r12})}{=} a_{j-1,k}^{j}(U)
    \stackrel{(\ref{eq-r8})}{=} a_{j-1,k-1}^{n}(U), \hspace{1em} i \leq j-1. \label{eq-r7}
\end{gather}
Hence it suffices to show that $a_{j,k}^{n}(U) = 0$ for $j < k-1$ and $k \leq n-1$.
The reason why we consider only $k \leq n-1$ is that an element type of $a_{j,n}^n(U)$ never appear in the above.
Then, observing the coefficient of $\bm{e}_{j+1, k}^n$ in (\ref{eq-r3}) for $l=j$, we obtain the required result.

\vspace{0.5em}

Next, we consider the other cases. If $j=k-1$, by the same argument as (\ref{eq-r7}) and (\ref{eq-r9}), we have
\[ a_{k-1,k}^{i}(U) = \begin{cases}
                              a_{k-1,k}^{k-2}(U) \stackrel{(\ref{eq-r7})}{=} a_{k-2,k-1}^{n}(U) = 0, \hspace{1em} & i \leq k-2, \\
                              a_{k-1,k}^{k+1}(U) \stackrel{(\ref{eq-r9})}{=} a_{k-1,k}^{n}(U) = 0, \hspace{1em} & i \geq k+1.
                           \end{cases}\]
For the case where $i= j$, $k$, we prepare some relations as follows.
By the coefficients of $\bm{e}_{j,l}^l$ and $\bm{e}_{l,k}^l$ in (\ref{eq-r3}), we see
\begin{gather}
 a_{j,l}^l(U) = a_{j,l+1}^{l+1}(U), \hspace{1em} 3 \leq j < l \leq n-1, \label{eq-r10} \\
 a_{l,k}^l(U) = a_{l+1,k}^{l+1}(U), \hspace{1em} 3 \leq l < k-1 \leq n-1 \label{eq-r11} 
\end{gather}
respectively. Then by (\ref{eq-r10}) and (\ref{eq-r11}),
\[ a_{j,i}^i(U) = a_{j,j+1}^{j+1}(U), \hspace{1em} a_{i,k}^i(U) = a_{k-1,k}^{k-1}(U). \]
Hence it suffices to show that $a_{l,l+1}^l(U) = a_{l,l+1}^{l+1}(U) = 0$ for $3 \leq l \leq n-1$.
From the coefficients of $\bm{e}_{l,l+1}^l$ in (\ref{eq-r3}) and (\ref{eq-r2}), we have
\[\begin{split}
     a_{l,l+1}^l(U) = -a_{l,l+1}^{l+1}(U), \hspace{1em}  a_{l,l+1}^{l+1}(U) = 0
  \end{split}\]
for $3 \leq l \leq n-1$ respectively. This shows the required result.

\vspace{1em}

{\bf (vi)} (Proof for $a_{2,k}^i(U) = 0$ for $3 \leq i, k \leq n$.)
First, we consider the case where $i \neq k$. For $k \neq l+1$,
in the equation (\ref{eq-r2}), the coefficient of $\bm{e}_{2,k}^l$ of $(Q^{-(l-1)} U Q^{l-1} - 1) f(U)$ is $a_{2,k}^{l+1}(U)$.
On the other hand, that of $(U-1) Q^{-(l-1)} f(U)$ is equal to the coefficient of $\bm{e}_{1,k}^l$ of $Q^{-(l-1)} f(U)$ times $-1$, and hence
$a_{k-l+1, n+2-l}^1(U)$. Therefore, we have
\begin{equation}
 a_{2,k}^{l+1}(U) = a_{k-l+1, n+2-l}^1(U) = 0 \label{eq-r200}
\end{equation}
for $3 \leq l \leq n-1$ and $k \neq l+1$.

\vspace{0.5em}

Similarly, from the coefficient of $\bm{e}_{2,l+1}^3$ in (\ref{eq-r2}), we see
\[ a_{2,l}^3(U) = - a_{2,n+2-l}^{n+4-l}(U) = 0 \]
from (\ref{eq-r200}) for $4 \leq l \leq n-1$. To show $a_{2,n}^3(U) = 0$, we consider the relation {\bf (N8)}:
\[ [U, Q^{-(n-2)}PUP^{-1}Q^{(n-2)}]=1. \]
Using (\ref{eq-comm}), we have
\[\begin{split}
   (Q^{-(n-2)} & PUP^{-1}Q^{(n-2)}-1)f(U) \\
      & = (U-1) f(Q^{-(n-2)} PUP^{-1}Q^{(n-2)}) \\
      & = (U-1)Q^{-(n-2)}(f(P)+Pf(U)-PUP^{-1}f(P)).
  \end{split}\]
Here the actions of $Q^{-(l-2)}PUP^{-1}Q^{(l-2)}$ for $2 \leq l \leq n$ on $e_k$ and $e_k^*$ are given by
\[\begin{split}
   Q^{-(l-2)}PUP^{-1}Q^{(l-2)} & \cdot e_k = \begin{cases}
                      e_l - e_{l-1}, \hspace{0.5em} & k=l, \\
                      e_k, \hspace{0.5em} & k \neq l,
                    \end{cases} \\
   Q^{-(l-2)}PUP^{-1}Q^{(l-2)} & \cdot e_k^*  = \begin{cases}
                      e_{l-1}^* + e_l^*, \hspace{0.5em} & k=l-1, \\
                      e_k^*, \hspace{0.5em} & k \neq l-1.
                    \end{cases}
  \end{split}\]
Then the coefficient of $\bm{e}_{2,n-1}^3$ of $(Q^{-(n-2)} PUP^{-1}Q^{(n-2)}-1)f(U)$ is given by $-a_{2,n}^3(U)$.
On the other hand, that of $(U-1)Q^{-(n-2)}(f(P)+Pf(U)-PUP^{-1}f(P))$ is equal to
\begin{center}
 the coefficient of $\bm{e}_{1,n-1}^3$ of $Q^{-(n-2)}(f(P)+Pf(U)-PUP^{-1}f(P))$ times $-1$,
\end{center}
and to
\begin{center}
 the coefficient of $\bm{e}_{2,3}^5$ of $f(P)+Pf(U)-PUP^{-1}f(P)$,
\end{center}
and to
\[ a_{2,3}^5(P) + a_{1,3}^5(U) - a_{2,3}^5(P) = a_{1,3}^5(U). \]
Hence we obtain $a_{2,n}^3(U) = - a_{1,3}^5(U) = 0$.

\vspace{0.5em}

Finally, we consider the case where $i=k$. In (\ref{eq-r2}), the coefficients of $\bm{e}_{2,l+1}^l$ of $(Q^{-(l-1)} U Q^{l-1} - 1) f(U)$
is equal to
\[ a_{2,l+1}^{l+1}(U) - a_{2,l}^{l}(U) - a_{2,l}^{l+1}(U) = a_{2,l+1}^{l+1}(U) - a_{2,l}^{l}(U). \]
On the other hand, that of $(U-1) Q^{-(l-1)} f(U)$ is equal to
\begin{center}
 the coefficient of $\bm{e}_{1,l+1}^l$ of $Q^{-(l-1)} f(U)$ times $-1$,
\end{center}
and to $a_{2, 2-l+n}^1(U) = 0$. Then we obtain $a_{2,l+1}^{l+1}(U) = a_{2,l}^{l}(U)$ for $3 \leq l \leq n-1$, and hence
\[ a_{2,n}^n(U) = a_{2,n-1}^{n-1}(U) = \cdots = a_{2,3}^3(U) = 0. \]

\vspace{0.5em}

Therefore we see $a_{j,k}^i(U) = 0$ for any $(i, j, k) \in I$. This shows that $f(U) = 0$. This completes the proof of Step III.

\vspace{1em}

\begin{flushleft}
{\bf Step IV.} (Proof for $a_{j,k}^i(P) = 0$.)
\end{flushleft}
By the results obtained in Step II, it suffices to show $a_{j,k}^i(P) = 0$ for $i \neq 2$.

\vspace{0.5em}

{\bf (i)} (Proof for $a_{j,k}^1(P) = 0$ for $1 \leq j < k \leq n$.)
From the relation {\bf (N12)}: $(PQ^{-1}UQ)^2=UQ^{-1}UQU^{-1}$ and a result $f(Q) = f(U) = 0$ as above, we see
\[ (1+PQ^{-1}UQ)f(PQ^{-1}UQ) = f(UQ^{-1}UQU^{-1}) = 0, \]
and hence
\begin{equation}\label{eq-q1}
 (1 + PQ^{-1}UQ)f(P) = 0.
\end{equation}
The actions of $PQ^{-1}UQ$ on $e_k$ and $e_k^*$ are given by
\[ PQ^{-1} UQ \cdot e_k = \begin{cases}
                      e_2, \hspace{0.5em} & k=1, \\
                      e_1 - e_3, \hspace{0.5em} & k=2, \\
                      e_k, \hspace{0.5em} & k \neq 1, 2,
                    \end{cases} \hspace{1em}
   PQ^{-1} UQ \cdot e_k^*  = \begin{cases}
                      e_2^*, \hspace{0.5em} & k=1, \\
                      e_1^*, \hspace{0.5em} & k=2, \\
                      e_1^* + e_3^*, \hspace{0.5em} & k=3, \\
                      e_k^*, \hspace{0.5em} & k \neq 1, 2, 3.
                    \end{cases} \]
Using this, we see that the coefficients of $\bm{e}_{2,3}^2$ and $\bm{e}_{3,k}^2$ in (\ref{eq-q1}) are calculated as
\begin{gather*}
 a_{2,3}^2(P) + a_{1,3}^1(P) - a_{1,2}^1(P) = 0, \\
 a_{3,k}^2(P) + a_{3,k}^1(P) - a_{2,k}^1(P) = 0, \hspace{1em} 4 \leq k \leq n
\end{gather*}
respectively.
Then from (\ref{eq-p2}) and (\ref{eq-p1}), we obtain $a_{1,2}^1(P) = 0$ and $a_{2,k}^1(P) = 0$ for
$4 \leq k \leq n$.

\vspace{0.5em}

Next, applying (\ref{eq-comm}) for $\sigma=P$ and $\tau=Q^{-(l-1)} U Q^{l-1}$ for $3 \leq l \leq n-1$, we have
\begin{equation}
  (Q^{-(l-1)} U Q^{l-1} - 1) f(P) = (U-1) f(Q^{-(l-1)} U Q^{l-1}) = 0. \label{eq-q2}
\end{equation}
By the coefficient of $\bm{e}_{2,l+1}^1$ in the equation above, we see $-a_{2,l}^1(P) = 0$.
In particular, $a_{2,3}^1(P) = 0$.
Furthermore, from the coefficient of $\bm{e}_{1,l+1}^1$, we see
\[ a_{1,l}^1(P) = 0, \hspace{1em} 3 \leq l \leq n-1. \]

\vspace{0.5em}

Now, from the relation {\bf (N2)}: $(QP)^{n-1}=1$,
\[ (1 + QP + \cdots + (QP)^{n-2})(f(Q) +Q f(P)) = 0, \]
and hence
\begin{equation}\label{eq-q3}
 (1 + QP + \cdots + (QP)^{n-2}) Q f(P) = 0.
\end{equation}
By the coefficient of $\bm{e}_{1,2}^1$ in (\ref{eq-q3}), we see
\[ a_{2,3}^2(P) +a_{2,4}^2(P) + \cdots + a_{2,n}^2(P) - a_{1,2}^2(P) = 0. \]
Then using (\ref{eq-p2}) and (\ref{eq-p3}), we obtain $a_{1,n}^1(P) = 0$.

\vspace{0.5em}

Subsequently, we consider $a_{j,k}^1(P)$ for $3 \leq j < k \leq n$.
From the coefficient of $\bm{e}_{j,l+1}^1$ in (\ref{eq-q2}) for $3 \leq j < l \leq n-1$, we obtain
\begin{equation}\label{eq-r19}
 a_{j,l}^1(P) = 0.
\end{equation}
On the other hand, from the coefficient of $\bm{e}_{l+1,n}^1$ in (\ref{eq-q2}), we see
$a_{l,n}^1(P) = 0$ for $3 \leq j \leq n-2$. Furthermore, observing the coefficient of $\bm{e}_{n-1,n}^1$ in (\ref{eq-q3}), we see
\[ -a_{1,2}^2(P) + a_{1,3}^2(P) +a_{3,4}^2(P) + \cdots + a_{n-2,n-1}^2(P) + a_{n-1,n}^2(P) = 0, \]
and hence $a_{n-1,n}^1(P) = - a_{n-1,n}^2(P) = 0$.
Therefore we have $a_{j,k}^1(P) = a_{j,k}^2(P) = 0$ for any $1 \leq j < k \leq n$.

\vspace{1em}

{\bf (ii)} (Proof for $a_{1,k}^i(P) = 0$ for $3 \leq i \leq n$ and $2 \leq  k \leq n$.)
First, we consider the case where $i \neq k$.
Observing the coefficients of $\bm{e}_{1,k}^l$ and $\bm{e}_{1,l+1}^3$ in (\ref{eq-q2}), we see
\begin{gather}
 a_{1,k}^{l+1}(P) = 0, \hspace{1em} 3 \leq l \leq n-1, \,\,\, 2 \leq k \neq l+1, \label{eq-q4} \\
 a_{1,l}^3(P) = 0, \hspace{1em} 4 \leq l \leq n-1 \label{eq-q6}
\end{gather}
respectively. Hence, it suffices to show that $a_{1,2}^3(P) = a_{1,n}^3(P) = 0$.

\vspace{0.5em}

By the coefficient of $\bm{e}_{1,2}^3$ in (\ref{eq-q3}), we see
\[ a_{2,3}^4(P) + a_{2,4}^5(P) + \cdots + a_{2,n-1}^n(P) + a_{2,n}^1(P) - a_{1,2}^3(P) = 0. \]
From (\ref{eq-p5}),
\[ - a_{1,3}^4(P) - a_{1,4}^5(P) - \cdots - a_{1,n-1}^n(P) + a_{2,n}^1(P) = a_{1,2}^3(P), \]
and hence $a_{1,2}^3(P) = 0$. Similarly, the coefficient of $\bm{e}_{1,n}^3$ in (\ref{eq-q3}), we see
\[ - a_{1,2}^4(P) + a_{2,3}^5(P) + \cdots + a_{2,n-2}^n(P) + a_{2,n-1}^1(P) + a_{1,n}^3(P) = 0, \]
and
\[ - a_{1,2}^4(P) - a_{1,3}^5(P) - \cdots - a_{1,n-2}^n(P) + a_{2,n-1}^1(P) = a_{1,n}^3(P), \]
and hence $a_{1,n}^3(P) = 0$.

\vspace{0.5em}

Next, we consider the case where $i=k$. By the coefficient of  $\bm{e}_{1,l+1}^l$ in (\ref{eq-q2}),
\begin{equation}
a_{1,l+1}^{l+1}(P) = a_{1,l}^l(P) - a_{1,l}^{l+1}(P) \stackrel{(\ref{eq-q4})}{=} a_{1,l}^l(P),
    \hspace{1em} 3 \leq l \leq n-1, \label{eq-q5} \\
\end{equation}
Hence it suffices to show that $a_{1,3}^3(P) = 0$.
On the other hand, by the coefficient of $\bm{e}_{2,3}^1$ in (\ref{eq-q1}), we see
\[ a_{2,3}^1(P) + a_{1,3}^2(P) - a_{1,2}^2(P) - a_{1,2}^3(P) + a_{1,3}^3(P) = 0, \]
and hence $a_{1,3}^3(P) = 0$.

\vspace{0.5em}

From the argument above, we obtain $a_{1,k}^i(P) = 0$ for $3 \leq i \leq n$ and $2 \leq  k \leq n$.
We remark that this also shows that $a_{2,k}^i(P) = 0$ for $3 \leq i \leq n$ and $3 \leq  k \leq n$ by (\ref{eq-p5}).

\vspace{1em}

{\bf (iii)} (Proof for $a_{j,k}^i(P) = 0$ for $3 \leq i \leq n$ and $3 \leq j <  k \leq n$.)
First, we consider the case where $i \geq 4$. By the coefficient of $\bm{e}_{j,k}^l$ in (\ref{eq-q2}), we see
\begin{equation}
a_{j,k}^{l+1}(P) = 0 \label{eq-q7}
\end{equation}
for $3 \leq l \leq n-1$ and $j, k \neq l+1$. Hence $a_{j,k}^i(P) = 0$ for $4 \leq i \leq n$ and $i \neq j, k$.

\vspace{0.5em}

If $i=j$ or $i=k$, observe the coefficients of $\bm{e}_{l,l+1}^l$, $\bm{e}_{j,l+1}^l$ and $\bm{e}_{l+1,k}^l$ in (\ref{eq-q2}).
Then we see
\begin{gather}
 a_{l,l+1}^{l+1}(P) = 0, \hspace{1em} 3 \leq l \leq n-1, \,\,\, \label{eq-q8-q9} \\
 a_{j,l+1}^{l+1}(P)-a_{j,l}^l(P) - a_{j,l}^{l+1}(P) = 0, \hspace{1em} 3 \leq l \leq n-1, \,\,\, j<l \label{eq-q8} \\
 a_{l+1,k}^{l+1}(P)-a_{l,k}^l(P) - a_{l,k}^{l+1}(P) = 0, \hspace{1em} 3 \leq l \leq n-1, \,\,\, l+1 < k \label{eq-q9}
\end{gather}
respectively. By (\ref{eq-q7}), the equations (\ref{eq-q8}) and (\ref{eq-q9}) are equivalent to
\[ a_{j,l+1}^{l+1}(P) = a_{j,l}^l(P), \hspace{1em} a_{l+1,k}^{l+1}(P) = a_{l,k}^l(P) \]
respectively. Using this and (\ref{eq-q8-q9}), we see
\[\begin{split}
   a_{j,n}^n(P) & = a_{j,n-1}^{n-1}(P) = \cdots = a_{j,j+1}^{j+1}(P) = 0\\
   a_{j,j+1}^j(P) & = a_{j-1,j+1}^{j-1}(P) = \cdots = a_{3,j+1}^3(P)
   \end{split}\]
for $3 \leq j \leq n-1$. Hence the proof of Step IV is finished if we show $a_{j,k}^3(P) = 0$ for $3 \leq j <  k \leq n$.

\vspace{0.5em}

From the coefficients of $\bm{e}_{j,l+1}^3$ and $\bm{e}_{l+1,n}^3$ in (\ref{eq-q2}), we see
\[\begin{split}
    a_{j,l}^3(P) & = 0, \hspace{1em} 3 \leq j < l \leq n-1, \\
    a_{l,n}^3(P) & = 0, \hspace{1em} 4 \leq l \leq n-2
  \end{split}\]
respectively. Hence it suffices to show $a_{3,n}^3(P) = a_{n-1,n}^3(P) = 0$.
Then observing the coefficients of $\bm{e}_{3,n}^3$ and $\bm{e}_{n-1,n}^3$ in (\ref{eq-q3}), we obtain
\[\begin{split}
   & -a_{1,4}^4(P) - a_{3,5}^5(P) - \cdots - a_{n-2,n}^n(P) + a_{1,n-1}^1(P) + a_{3,n}^3(P) = 0, \\
   & -a_{1,n}^4(P) + a_{1,3}^5(P) + a_{3,4}^6(P) + \cdots + a_{n-3,n-2}^n(P)+a_{n-2,n-1}^1(P) +a_{n-1,n}^3(P) = 0.
  \end{split}\]
These equations induce the required results. Therefore we obtain $a_{j,k}^i(P) = 0$ for any $(i, j, k) \in I$.
This shows that $f(P) = 0$. This completes the proof of Step IV.

\vspace{0.5em}

\begin{flushleft}
{\bf Step V.} (The rest of the proof for $a_{j,k}^i(S) = 0$.)
\end{flushleft}
Here we show that $a_{1,k}^i(S) = 0$ for $i, k \geq 2$, and $a_{j,k}^1(S) = 0$ for $2 \leq j < k \leq n$.

\vspace{0.5em}

By the relation {\bf (N11)}: $SUSPS=PU^{-1}PU$, we have
\begin{equation}\label{eq-s1}
 (1+SU+SUSP)f(S) = f(PU^{-1}PU) = 0.
\end{equation}
The actions of $SU$ and $SUSP$ on $e_k$ and $e_k^*$ are given by
\[ SU \cdot e_k = \begin{cases}
                      - e_1- e_2, \hspace{0.5em} & k=1, \\
                      e_k, \hspace{0.5em} & k \neq 1,
                    \end{cases} \hspace{1em}
   SU \cdot e_k^*  = \begin{cases}
                      - e_1^*, \hspace{0.5em} & k=1, \\
                      - e_1^* +e_2^*, \hspace{0.5em} & k=2, \\
                      e_k^*, \hspace{0.5em} & k \neq 1, 2,
                    \end{cases} \]
\[ SUSP \cdot e_k = \begin{cases}
                      e_2, \hspace{0.5em} & k=1, \\
                      e_1 + e_2, \hspace{1em} & k=2, \\
                      e_k, \hspace{0.5em} & k \neq 1, 2
                    \end{cases} \hspace{1em}
   SUSP \cdot e_k^*  = \begin{cases}
                      - e_1^* + e_2^*, \hspace{0.5em} & k=1, \\
                      e_1^*, \hspace{0.5em} & k=2, \\
                      e_k^*, \hspace{0.5em} & k \neq 1, 2,
                    \end{cases} \]
Using this, for $3 \leq k \leq n$, from the coefficients $\bm{e}_{1,k}^1$ and $\bm{e}_{1,k}^2$ in (\ref{eq-s1}), we obtain
\[\begin{split}
   2 a_{1,k}^1(S) + a_{1,k}^2(S) & + a_{2,k}^2(S) - a_{2,k}^1(S) = 0, \\
   a_{2,k}^1(S) & = 0 \\
\end{split}\]
respectively. Hence $a_{1,k}^2(S) = a_{2,k}^1(S) = 0$.

\vspace{0.5em}

Next, we show $a_{1,k}^i(S) = 0$ for $i \geq 3$ and $k \geq 2$.
Consider the relator {\bf (N4)}, and applying (\ref{eq-comm}) for $\sigma=S$ and $\tau=QP$, we have
\begin{equation}
    (QP - 1) f(S) = (S-1) f(QP) = 0. \label{eq-qp1}
\end{equation}
The actions of $QP$ on $e_k$ and $e_k^*$ are given by
\[ QP \cdot e_k = \begin{cases}
                      e_1, \hspace{0.5em} & k=1, \\
                      e_n, \hspace{0.5em} & k=2, \\
                      e_{k-1}, \hspace{0.5em} & k \neq 1, 2,
                    \end{cases} \hspace{1em}
   QP \cdot e_k^*  = \begin{cases}
                      e_1^*, \hspace{0.5em} & k=1, \\
                      e_n^*, \hspace{0.5em} & k=2, \\
                      e_{k-1}^*, \hspace{0.5em} & k \neq 1, 2.
                    \end{cases} \]
Using this, from the coefficients $\bm{e}_{1,k}^i$ in (\ref{eq-qp1}), we see
\begin{gather}
  a_{1,k+1}^{i+1}(S) = a_{1,k}^i(S), \hspace{1em} 2 \leq i, k \leq n-1, \label{eq-s3} \\
  a_{1,2}^{i+1}(S) = a_{1,n}^i(S), \hspace{1em} 2 \leq i \leq n-1, \,\,\, k=n. \label{eq-s4}
\end{gather}
From (\ref{eq-s3}), if $i \leq k$,
\[ a_{1,k}^i(S) = a_{1,k-1}^{i-1}(S) = \cdots = a_{1,k+2-i}^{2}(S) = 0, \]
and hence from (\ref{eq-s3}) if $i > k$,
\[ a_{1,k}^i(S) = a_{1,k-1}^{i-1}(S) = \cdots = a_{1,2}^{i+2-k}(S) = 0.\]

\vspace{0.5em}

Finally, we show $a_{j,k}^1(S) = 0$ for $3 \leq j < k \leq n$.
From the coefficient $\bm{e}_{j,k}^1$ in (\ref{eq-qp1}), we see
\[ a_{j+1,k+1}^{1}(S) = a_{j,k}^1(S) \]
for $2 \leq j < k \leq n-1$. This shows that for $3 \leq j < k \leq n$,
\[ a_{j,k}^1(S) = a_{j-1,k-1}^1(S) = \cdots = a_{2,k+2-j}^1(S) = 0. \]

\vspace{0.5em}

Therefore we obtain $a_{j,k}^i(S) = 0$ for any $(i, j, k) \in I$. This completes the proof of Step V.

\vspace{0.5em}

From the argument above, we verify that the map $\Phi : \mathrm{Cros}(\mathrm{Aut}\,F_n, V_L) \rightarrow L^{\oplus (n^3-n^2+4)/2}$
is injective. Through the map $\Phi$, consider the $L$-modlue $\mathrm{Cros}(\mathrm{Aut}\,F_n, V_L)$
as a submodule of $W:=L^{\oplus (n^3-n^2+4)/2}$.
Next, we study the quotient $L$-module $W/\mathrm{Prin}(\mathrm{Aut}\,F_n, V_L)$.

\vspace{0.5em}

\noindent
{\bf Claim 2.} $W/\mathrm{Prin}(\mathrm{Aut}\,F_n, V_L)$ is a free $L$-module of rank $2$.

\vspace{0.5em}

For any element
\[ \bm{a} := \sum_{(i,j,k) \in I} a_{j,k}^i \bm{e}_{j,k}^i \in V_L, \]
let $f_{\bm{a}} : \mathrm{Aut}\,F_n \rightarrow V_L$ be the principal crossed homomorphisms associated to $\bm{a}$. 
For example,
\[\begin{split}
   f_{\bm{a}}(Q) & = \sum_{i \neq n,\, 1 \leq j < k \leq n-1} (a_{j+1, k+1}^{i+1}-a_{j,k}^i) \bm{e}_{j,k}^i 
              + \sum_{1 \leq i, j \leq n-1} (-a_{1, j+1}^{i+1}-a_{j,n}^i) \bm{e}_{j,n}^i \\
          & \hspace{1.5em} + \sum_{1 \leq j < k \leq n} (a_{j+1,k+1}^1-a_{j,k}^n) \bm{e}_{j,k}^n
              + \sum_{1 \leq j \leq n-1} (-a_{1, j+1}^{1}-a_{j,n}^n) \bm{e}_{j,n}^n \\
  \end{split}\]
and
\[\begin{split}
    f_{\bm{a}}(U) & = \sum_{2 \leq k \leq n} a_{1,k}^2 \bm{e}_{1,k}^1 + \sum_{3 \leq k \leq n} (a_{2,k}^2-a_{1,k}^1-a_{1,k}^2) \bm{e}_{2,k}^1
              + \sum_{3 \leq j < k \leq n} a_{j,k}^2 \bm{e}_{j,k}^1  \\
          & \hspace{1.5em} + \sum_{3 \leq k \leq n} - a_{1,k}^2 \bm{e}_{2,k}^2
              + \sum_{3 \leq i, k \leq n} -a_{1,k}^i \bm{e}_{2,k}^i.
  \end{split}\]
In order to determine the $L$-module structure of $W/\mathrm{Prin}(\mathrm{Aut}\,F_n, V_L)$, it
suffices to find the elementary divisors of an $n^2(n-1)/2 \times (n^3-n^2+4)/2$ matrix:
\[\begin{split} A := \bordermatrix{
             &  a_{j,k}^1(Q) & a_{j,k}^2(Q) & \cdots & a_{j,k}^{n-1}(Q) & a_{j,k}^1(U) & a_{1,2}^2(S) & a_{2,3}^3(U) \\
   a_{j,k}^1 &  A^{1,1}      & A^{1,2}      & \cdots & A^{1,n-1}        & A^{1,n}      & A^{1,n+1}    & A^{1,n+2}    \\
   a_{j,k}^2 &  A^{2,1}      & A^{2,2}      & \cdots & A^{2,n-1}        & A^{2,n}      & A^{2,n+1}    & A^{2,n+2}    \\
 \,\,\, \vdots &  \vdots       & \vdots       & \vdots & \vdots           & \vdots       & \vdots       & \vdots       \\
   a_{j,k}^n &  A^{n,1}      & A^{n,2}      & \cdots & A^{n,n-1}        & A^{n,n}      & A^{n,n+1}    & A^{n,n+2}    \\
                                   }  \end{split}\]
whose row is indexed by $a_{j,k}^i$s, and whose column is indexed by (\ref{eq-gen}).
Here each $A^{p,q}$ is a block matrix defined as follows. First, we consider the case where $1 \leq q \leq n-1$.
For $\sigma =P, Q, S$ and $U$, set
\[ f_{\bm{a}}(\sigma) := \sum_{(i, j, k) \in I} a_{j,k}^i(\sigma) \bm{e}_{j,k}^i \in V_L. \]
Then, for any $1 \leq j_2 < k_2 \leq n$, we have
\[ a_{j_2, k_2}^q(Q) = \sum_{(p, j_1, k_1) \in I} C_{(j_1,k_1), \,\, (j_2,k_2)}^{p,q} a_{j_1,k_1}^p \]
for some $C_{(j_1,k_1), \,\, (j_2,k_2)}^{p,q} \in L$. Then the matrix $A^{p,q}$ is defined by
\[\begin{split} A^{p,q} := 
      \bordermatrix{
             & a_{1,2}^q(Q)                  & a_{1,3}^q(Q)                & \cdots      & a_{n-1,n}^q(Q)                   \\
   a_{1,2}^p & C_{(1,2), \,\, (1,2)}^{p,q}   & C_{(1,2), \,\, (1,3)}^{p,q} & \cdots      & C_{(1,2), \,\, (n-1,n)}^{p,q}    \\
   a_{1,3}^p & C_{(1,3), \,\, (1,2)}^{p,q}   & C_{(1,3), \,\, (1,3)}^{p,q} & \cdots      & C_{(1,3), \,\, (n-1,n)}^{p,q}    \\
\,\,\, \vdots & \vdots                        & \vdots                      & \vdots      & \vdots                           \\
 a_{n-1,n}^p & C_{(n-1,n), \,\, (1,2)}^{p,q} & C_{(n-1,n), \,\, (1,3)}^{p,q} & \cdots    & C_{(n-1,n), \,\, (n-1,n)}^{p,q}  \\
                   }  \end{split}\]
where the rows are indexed by $a_{j,k}^p$s according to the usual lexicographic order on the set
$\{ (j,k) \,|\, 1 \leq j < k \leq  n \}$. Similarly, the columns are indexed by $a_{j,k}^q(Q)$s.

\vspace{0.5em}

By an argument similar to the above,
the block matrices $A^{p,n}$, $A^{p,n+1}$ and $A^{p,n+2}$ for $1 \leq p \leq n$ are defined from
$a_{j,k}^1(U)$s, $a_{1,2}^2(S)$ and $a_{2,3}^3(U)$ respectively.

\vspace{0.5em}

Set
{\small
\[\begin{split} A' :=  \bordermatrix{
             &  a_{j,k}^1(Q) & a_{j,k}^2(Q) & \cdots & a_{j,k}^{n-1}(Q) & a_{j,k}^1(U) \\
   a_{j,k}^1 &  A^{1,1}      & A^{1,2}      & \cdots & A^{1,n-1}        & A^{1,n}      \\
   a_{j,k}^2 &  A^{2,1}      & A^{2,2}      & \cdots & A^{2,n-1}        & A^{2,n}      \\
\,\,\, \vdots  &  \vdots       & \vdots       & \vdots & \vdots           & \vdots       \\
   a_{j,k}^n &  A^{n,1}      & A^{n,2}      & \cdots & A^{n,n-1}        & A^{n,n}      \\
               }.  \end{split}\]
}
In the following, we prove that all elementary divisors of $A$ are equal to $1 \in L$ by showing that $A'$ can be transformed into the
identity matrix with only the elementary column operations. Then
we conclude $W/\mathrm{Prin}(\mathrm{Aut}\,F_n, V_L) \cong L^{\oplus 2}$.

\vspace{0.5em}

First, we transform the $a_{j,k}^1(U)$ columns of $A'$. To do this, we use the followings. From $f_{\bm{a}}(Q)$ and $f_{\bm{a}}(U)$ as above, 
the $a_{1,k}^1(Q)$ columns,
the $a_{j,k}^1(U)$ columns for $j \neq 2$ and the $a_{2,k}^1(U)$ columns of $A'$ are given by
{\small
\[\begin{split}
\bordermatrix{
                 &  a_{1,k}^1(Q) & a_{1,n}^1(Q)  \\
     a_{1,k}^1   & -E            & O             \\
     a_{1,n}^1   & O             & -1            \\
     a_{1,2}^2   & O             & -1            \\
     a_{1,k+1}^2   & O             & O             \\
     a_{2,k+1}^2 & E             & O             \\
\,\,\, \vdots    & O             & O             \\
          }, \hspace{1em}
\bordermatrix{
                 &  a_{j,k}^1(U) \\
\,\,\, \vdots    & O             \\
     a_{j,k}^{2} & E             \\
\,\,\, \vdots    & O             \\
          }, \hspace{1em}
\bordermatrix{
                  &  a_{2,k}^1(U) \\
     a_{1,2}^1    &  O            \\
     a_{1,k}^1    & -E            \\
\,\,\, \vdots     & O             \\
     a_{1,k}^2    & -E            \\
     a_{2,k}^2    & E             \\
\,\,\, \vdots     & O             \\
          } \end{split}\]
}
respectively. Here $E$ denotes the identity matrix.

\vspace{0.5em}

Let us consider the $a_{2,k}^1(U)$ columns of $A'$. Add the $a_{1,k}^1(U)$ columns
to $a_{2,k}^1(U)$ columns for $3 \leq k \leq n$, and minus $a_{1,k}^1(Q)$ columns from the $a_{2,k}^1(U)$ columns for $3 \leq k \leq n$.
Subsequently, by subtracting the $a_{1,2}^1(U)$ column from the $a_{2,n}^1(U)$ column, we see that the $a_{2,k}^1(U)$ columns of $A'$ are transformed into
{\small
\[\begin{split}
      \bordermatrix{
             &   a_{2,k}^1(U) \\
\,\,\, \vdots&   O            \\
   a_{2,k}^2 &   X            \\
\,\,\, \vdots &  O            \\
                   } \hspace{1em} \text{where} \hspace{1em}
    X = \bordermatrix{
             & a_{2,3}^1(U)  & a_{2,4}^1(U)  & \cdots    & a_{2,n-1}^1(U)    & a_{2,n}^1(U)  \\
   a_{2,3}^2 & 1             & 0             & \cdots    & 0                 & 0               \\
   a_{2,4}^2 & -1            & 1             &           &                   & 0              \\
\,\,\, \vdots  & 0             & -1            & \ddots    &                   & \vdots          \\
\,\,\, \vdots  & \vdots        &               & \ddots    & 1                 & 0              \\
   a_{2,n}^2 & 0             & 0             &           & -1                & 1               \\
            }.  \end{split}\]
}
It is easily seen that $X$ can be transformed into the identity matrix with the elementary column operations.
Hence the $a_{j,k}^1(U)$ columns of $A'$ are transformed into
\begin{equation}\label{eq-mat}\begin{split}
\bordermatrix{
               &  a_{j,k}^{1}(U) \\
\,\,\, \vdots  & O             \\
    a_{j,k}^2  & E             \\
\,\,\, \vdots  & O             \\
          }. \end{split}\end{equation}

\vspace{0.5em}

Next, for any $1 \leq j \leq n-1$, we consider the $a_{j,k}^1(Q)$ columns given by
{\small
\[\begin{split}
\bordermatrix{
                 &  a_{j,k}^1(Q) & a_{j,n}^1(Q)  \\
     a_{j,k}^1   & -E            & O             \\
     a_{j,n}^1   & O             & -1            \\
\,\,\, \vdots    & O             & O             \\
     a_{1,j+1}^2   & O             & -1            \\
\,\,\, \vdots    & O             & O             \\
     a_{j+1,k+1}^2 & E             & O             \\
\,\,\, \vdots    & O             & O             \\
          } \end{split}\]
}
Multiplying each column by $-1$ and using (\ref{eq-mat}), we can transform the $a_{j,k}^1(Q)$ columns into
\[\begin{split}
      \bordermatrix{
               &  a_{j,k}^1(Q) \\
     a_{j,k}^1 & E             \\
 \,\,\, \vdots & O             \\
          }. \end{split}\]

\vspace{0.5em}

Now, we consider the $a_{j,k}^2(Q)$ columns given by
{\scriptsize
\[\begin{split}
\bordermatrix{
                 &  a_{1,k}^2(Q) & a_{1,n}^2(Q) &  a_{2,l}^2(Q) & a_{2,n}^2(Q) &  \cdots       & a_{n-2,n-1}^2(Q) & a_{n-1,n}^2(Q) \\
\,\,\, \vdots    & O             & O            & O             & O            &  \cdots       & O              & O            \\
     a_{1,k}^2   & -E            & O            & O             & O            &  \cdots       & O              & O            \\
     a_{1,n}^2   & O             & -1           & O             & 0            &  \cdots       & 0              & 0            \\
     a_{2,l}^2   & O             & O            & -E            & O            &  \cdots       & O              & O            \\
     a_{2,n}^2   & O             & 0            & O             & -1           &  \cdots       & 0              & 0            \\
\,\,\, \vdots    & \vdots        & \vdots       & \vdots        & \vdots       &  \vdots       & \vdots         & \vdots       \\
   a_{n-1,n}^2   & O             & 0            & O             & 0            &  \cdots       & 0              & -1           \\
     a_{1,2}^3   & O             & -1           & O             & 0            &  \cdots       & 0              & 0            \\
     a_{1,3}^3   & O             & 0            & O             & -1           &  \cdots       & 0              & 0            \\
\,\,\, \vdots    & \vdots        & \vdots       & \vdots        & \vdots       &  \vdots       & \vdots         & \vdots       \\
   a_{1,n}^3     & O             & O            & O             & O            &  \cdots       & 0              & -1           \\
   a_{2,k+1}^3   & E             & O            & O             & 0            &  \cdots       & O              & O            \\
   a_{3,l+1}^3   & O             & O            & E             & O            &  \cdots       & O              & O            \\
\,\,\, \vdots    & \vdots        & \vdots       & \vdots        & \vdots       &  \vdots       & \vdots         & \vdots       \\
   a_{n-1,n}^3    & O             & O            & O             & O            &  \cdots       & 1              & 0            \\
\,\,\, \vdots    & O             & O            & O             & O            &  \cdots       & O              & O            \\
          }.
\end{split}\]
}
respectively. Similarly, using (\ref{eq-mat}) and the multiplication of $-1$, we can transform
the $a_{j,k}^2(Q)$ columns of $A'$ into
\[\begin{split}
      \bordermatrix{
               &  a_{j,k}^2(Q) \\
 \,\,\, \vdots & O             \\
     a_{j,k}^3 & E             \\
 \,\,\, \vdots & O             \\
          }. \end{split}\]
respectively. 

\vspace{0.5em}

By the same argument as above, for any $3 \leq p \leq n-1$, we can transform the $a_{j,k}^p(Q)$ columns of $A'$ into
\[\begin{split}
\bordermatrix{
               &  a_{j,k}^p(Q) \\
 \,\,\, \vdots & O             \\
     a_{j,k}^{p+1} & E             \\
 \,\,\, \vdots & O             \\
          } \end{split}\]
recursively from $p=3$ to $n-1$. From the argument above, we can transform $A'$ into 
{\small
\[\begin{split} A' :=  \bordermatrix{
                &  a_{j,k}^1(Q) & a_{j,k}^2(Q) & \cdots & a_{j,k}^{n-1}(Q) & a_{j,k}^1(U) \\
      a_{j,k}^1 &  E            & O            & \cdots & O                & O      \\
      a_{j,k}^2 &  O            & O            & \cdots & O                & E      \\
      a_{j,k}^3 &  O            & E            & \cdots & O                & O      \\
 \,\,\, \vdots  &  \vdots       & \vdots       & \ddots & \vdots           & \vdots       \\
      a_{j,k}^n &  O            & O            & \cdots & E                & O      \\
               },  \end{split}\]
}
and into the identity matrix with only the elementary column operations. Therefore we conclude that all elementary divisors of $A$ are equal to $1 \in L$.
In particular, observing the process of the transformation of $A$ as mentioned above, we see that a map $\Phi' : W \rightarrow L^{\oplus 2}$ defined by
\[ \Big{(} (a_{j,k}^i(Q))_{i \neq n, \,\, 1 \leq j < k \leq n}, \,\, (a_{j,k}^1(U))_{1 \leq j < k \leq n}, \,\, a_{1,2}^2(S), \,\, a_{2,3}^3(U) \Big{)}
   \mapsto (a_{1,2}^2(S), \,\, a_{2,3}^3(U)) \]
induces an isomorphism
\[ W/\mathrm{Prin}(\mathrm{Aut}\,F_n, V_L) \cong L^{\oplus 2}. \]

\vspace{1em}

Finally, for the crossed homomorphisms $f_M$, $f_N \in \mathrm{Cros}(F, V_L)$ defined in Section {\rmfamily \ref{S-Cro}}, we see
\[ \Phi'(f_M) =(-1,0), \hspace{1em} \Phi(f_N) = (0,1). \]
Hence
\[ H^1(\mathrm{Aut}\,F_n,V_L) \cong W/\mathrm{Prin}(\mathrm{Aut}\,F_n, V_L) \cong L^{\oplus 2}. \]
This completes the proof of Theorem {\rmfamily \ref{T-coh}}. $\square$

\vspace{1em}

From Theorem {\rmfamily \ref{T-coh}}, we see that the crossed homomorphisms $f_M$ and $f_N$, and hence $f_M$ and $f_K$,
generate $H^1(\mathrm{Aut}\,F_n,V_L)$ for $n \geq 5$.

\section{Some Applications}\label{S-App}

In this section, we consider the first Johnson homomorphism and the first cohomology group of the outer automorphism group.

\subsection{The first Johnson homomorphism}
\hspace*{\fill}\ 

\vspace{0.5em}

As a corollary to Theorem {\rmfamily \ref{T-coh}}, here we show
\begin{cor}\label{p-joh}
Let $L$ be a principal ideal domain without $2$-torsion. If $L$ does not contain $1/2$, then for $n \geq 5$,
there is no crossed homomorphism from $\mathrm{Aut}\,F_n$ to $V_L$ which restriction to
$\mathrm{IA}_n$ coincides with $\tau_{1,L}$.
\end{cor}

\textit{Proof of Corollary {\rmfamily \ref{p-joh}}.}
Assume that the first Johnson homomorphism $\tau_{1,L}$ extends to $\mathrm{Aut}\,F_n$ as a crossed homomorphism.
By Theorem {\rmfamily \ref{T-coh}}, $\tau_{1,L}$ is cohomologous to $a f_M +b f_K$ for some $a$, $b \in L$.
Then, for distinct $i$, $j$, $k$ and $j<k$, observing $f_K(K_{ijk}) = 2 \bm{e}_{j,k}^i$, we see
\[ \bm{e}_{j,k}^i = \tau_{1,L}(K_{ijk})= 2 b \bm{e}_{j,k}^i. \]
This shows $2b=1$. It is contradiction to the hypothesis of $L$. This completes the proof of Corollary {\rmfamily \ref{p-joh}}. $\square$

\vspace{0.5em}

\subsection{Outer automorphism group}\label{S-Out}
\hspace*{\fill}\ 

\vspace{0.5em}

Here, we compute the first cohomology group of the outer automorphism group
$\mathrm{Out}\,F_n := \mathrm{Aut}\,F_n/\mathrm{Inn}\,F_n$ of $F_n$ with coefficients in $V_L$
for any principal ideal domain $L$ without $2$-torsion.

\begin{pro}\label{P-out}
Let $L$ be as above. Then, for $n \geq 5$,
\[ H^1(\mathrm{Out}\,F_n, V_L) = L. \]
\end{pro}

\textit{Proof of Proposition {\rmfamily \ref{P-out}}.}
Considering the five-term exact sequence of
\[ 1 \rightarrow \mathrm{Inn}\,F_n \rightarrow \mathrm{Aut}\,F_n \rightarrow \mathrm{Out}\,F_n \rightarrow 1, \]
we have
\[ 0 \rightarrow H^1(\mathrm{Out}\,F_n,V_L) \xrightarrow{} H^1(\mathrm{Aut}\,F_n,V_L) \xrightarrow{\alpha} H^1(\mathrm{Inn}\,F_n,V_L)^{\mathrm{Out}\,F_n}. \]
From Proposition {\rmfamily \ref{T-coh}}, $H^1(\mathrm{Aut}\,F_n,V_L)$ can be identified with a free $L$-module $L^{\oplus 2}$ generated by $f_M$ and $f_K$.
Since
\[\begin{split}
   \alpha(f_M) & = (n-1) \sum_{i=1}^n \iota_i^* \otimes (\bm{e}_{1,i}^1 + \bm{e}_{2,i}^2 + \stackrel{\check{i}}{\cdots} + \bm{e}_{n,i}^n), \\
   \alpha(f_K) & = 2 \sum_{i=1}^n \iota_i^* \otimes (\bm{e}_{1,i}^1 + \bm{e}_{2,i}^2 + \stackrel{\check{i}}{\cdots} + \bm{e}_{n,i}^n),
  \end{split}\]
the image of $\alpha$ is contained in a free $L$-module $L$ generated by
\[ \sum_{i=1}^n \iota_i^* \otimes (\bm{e}_{1,i}^1 + \bm{e}_{2,i}^2 + \stackrel{\check{i}}{\cdots} + \bm{e}_{n,i}^n). \]
Then $\alpha$ is considered as an $L$-linear homomorphism $L^{\oplus 2} \rightarrow L$ which matrix representation is
$\begin{pmatrix} n-1 & 2 \end{pmatrix}$. Using the elementary operations, we can transform it into
\[\begin{split}
    & \begin{pmatrix} 1 & 0 \end{pmatrix}, \hspace{1em} n: \, \text{even}, \\
    & \begin{pmatrix} 2 & 0 \end{pmatrix}, \hspace{1em} n: \, \text{odd}.
  \end{split}\]
In both cases, the kernel of this homomorphism is isomorphic to $L$. This completes the proof of Proposition {\rmfamily \ref{P-out}}. 

\section{Acknowledgments}\label{S-Ack}

This research is supported by a JSPS Research Fellowship for Young Scientists and the Global COE program at Kyoto University.


\begin{thebibliography}{99}
 \bibitem{Bir} J. S. Birman; Braids, Links, and Mapping Class Groups, Annals of Math. Studies 82
               (Princeton University Press, 1974).
 \bibitem{Bou} N. Bourbaki; Lie groups and Lie algebra, Chapters 1--3, Softcover edition of the 2nd printing, Springer-Verlag (1989).
 \bibitem{Co1} F. Cohen and J. Pakianathan; On Automorphism Groups of Free Groups, and Their Nilpotent Quotients, preprint.
 \bibitem{Co2} F. Cohen and J. Pakianathan; On subgroups of the automorphism group of a free group and associated graded Lie algebras,
               preprint.
 \bibitem{Far} B. Farb; Automorphisms of $F_n$ which act trivially on homology, in preparation.
 \bibitem{Gal} S. Galatius; Stable homology of automorphism groups of free groups, preprint, \\ {\texttt{arXiv:math.AT/0610216v3}}.
 \bibitem{Ger} S. M. Gersten; A presentation for the special automorphism group of a free group,
               J. Pure and Applied Algebra 33 (1984), 269-279.
 \bibitem{HaV} A. Hatcher and K. Vogtmann; Rational homology of $\mathrm{Aut}(F_n)$, Math. Res.
               Lett. 5 (1998), 759-780.
 \bibitem{HaN} A. Hatcher and N. Wahl; Stabilization for the automorphisms of free groups with boundaries,
               Geometry and Topology, Vol. 9 (2005), 1295-1336.
 \bibitem{Kaw} N. Kawazumi; Cohomological aspects of Magnus expansions, preprint,
               {\texttt{arXiv:math.GT/0505497}}.
 \bibitem{Krs} S. Krsti\'{c}, J. McCool; The non-finite presentability in $IA(F_3)$ and $GL_{2}(\Z[t,t^{-1}])$, Invent. Math. 129
               (1997), 595-606.
 \bibitem{Mag} W. Magnus; $\ddot{\mathrm{U}}$ber $n$-dimensinale Gittertransformationen,
               Acta Math. 64 (1935), 353-367.
 \bibitem{MKS} W. Magnus, A. Karras and D. Solitar; Combinatorial group theory, Interscience Publ., New York (1966).
 \bibitem{McC} J. McCool; Some remarks on IA automorphisms of free groups, Can. J. Math. Vol. XL, no. 5 (1998), 1144-1155.
 \bibitem{Mo6} S. Morita; On the Homology Groups of the Mapping Class Groups of Orientable Surfaces with Twisted Coefficients,
               Proc. Japan Acad., 62, Ser. A (1986), 148-151.
 \bibitem{Mo0} S. Morita; Families of Jacobian manifolds and characteristic classes of surface
               bundles I, Ann. Inst. Fourier 39 (1989), 777-810.
 \bibitem{Mo4} S. Morita; Families of Jacobian manifolds and characteristic classes of surface bundles, II, Math. Proc. Camb. Phil. Soc. 105 (1989), 79-101.
 \bibitem{Mo1} S. Morita; Abelian quotients of subgroups of the mapping class group of surfaces, Duke Mathematical
               Journal 70 (1993), 699-726.
 \bibitem{Mo5} S. Morita; The extension of Johnson's homomorphism from the Torelli group to the mapping class group, Invent. math. 111 (1993), 197-224.
 \bibitem{Mo2} S. Morita; Structure of the mapping class groups of surfaces: a survey and a prospect, Geometry and Topology
               Monographs Vol. 2 (1999), 349-406.
 \bibitem{Mo3} S. Morita; Cohomological structure of the mapping class group and beyond, preprint.
 \bibitem{Ni0} J. Nielsen; Die Isomorphismen der allgemeinen unendlichen Gruppe mit zwei Erzeugenden, Math. Ann. 78 (1918),
               385-397.
 \bibitem{Ni1} J. Nielsen; Die Isomorphismengruppe der freien Gruppen, Math. Ann. 91 (1924), 169-209. 
 \bibitem{Ni2} J. Nielsen; Untersuchungen zur Topologie der geschlossenen Zweiseitigen Fl$\ddot{\mathrm{a}}$schen, Acta Math.
               50 (1927), 189-358.
 \bibitem{Sa1} T. Satoh; Twisted first homology group of the automorphism group of a free group, Journal of Pure and Applied Algebra,
               204 (2006), 334-348.
 \bibitem{Sat} T. Satoh; New obstructions for the surjectivity of the Johnson homomorphism of
               the automorphism group of a free group, Journal of the London Mathematical Society, (2) 74 (2006) 341-360.
 \bibitem{Sa2} T. Satoh; Twisted second homology group of the automorphism group of a free group, Journal of Pure and Applied Algebra,
               211 (2007), 547-565.
\end{thebibliography}
\end{document}